\documentclass[12pt,a4paper,final, dvipsnames]{amsart}

\input xy
\xyoption{all}

\DeclareMathAlphabet{\mathpzc}{OT1}{pzc}{m}{it}

\usepackage{amsfonts,amsmath,amssymb,indentfirst,mathrsfs,amscd}
\usepackage[mathscr]{eucal}
\usepackage{graphicx}
\usepackage{nicefrac}
\usepackage{amsthm}

\usepackage{color}
\usepackage[utf8]{inputenc}
\usepackage[spanish, english]{babel}
\usepackage{booktabs}
\usepackage{multirow}
\usepackage{verbatim}
\usepackage{hyperref}

\usepackage{calrsfs}
\DeclareMathAlphabet{ \mathcal}{OMS}{zplm}{m}{n}

\usepackage[all]{xy}
\usepackage{graphicx}
\usepackage{stmaryrd}
\usepackage{enumitem}

\usepackage{empheq}
\usepackage{color}
\usepackage{xcolor}

\author[]{\'Angel Gonz\'alez-Prieto}
\address{Instituto de Ciencias Matem\'aticas CSIC-UAM-UC3M-UCM. C.\ Nicol\'as Cabrera, 13-15, 28049. Madrid, Spain.}
\address{ETSI de Sistemas Inform\'aticos, Universidad Polit\'ecnica de Madrid. C.\ Alan Turing s/n, 28031. Madrid, Spain.}
\email{angel.gonzalez.prieto@icmat.es}
\title[]{Virtual classes of\\ parabolic $\SL{2}(\CC)$-character varieties}


\DeclareMathOperator{\coker}{coker\,}
\DeclareMathOperator{\tr}{Tr\,}             

\DeclareMathOperator{\Gr}{Gr}

\usepackage{amsmath}
\begin{document}

\newtheorem{thm}{Theorem}[section]
\newtheorem{prop}[thm]{Proposition}
\newtheorem{lem}[thm]{Lemma}
\newtheorem{cor}[thm]{Corollary}
\newtheorem{conjecture}{Conjecture}
\newtheorem*{theorem*}{Theorem}

\theoremstyle{definition}
\newtheorem{defn}[thm]{Definition}
\newtheorem{ex}[thm]{Example}
\newtheorem{as}{Assumption}

\theoremstyle{remark}
\newtheorem{rmk}[thm]{Remark}

\theoremstyle{remark}
\newtheorem*{prf}{Proof}

\newcommand{\iacute}{\'{\i}} 
\newcommand{\norm}[1]{\lVert#1\rVert} 

\newcommand{\lto}{\longrightarrow}
\newcommand{\hra}{\hookrightarrow}

\newcommand{\suchthat}{\;\;|\;\;}
\newcommand{\dbar}{\overline{\partial}}

\newcommand{\cA}{\mathcal{A}}
\newcommand{\cC}{\mathcal{C}}
\newcommand{\cD}{\mathcal{D}}
\newcommand{\cE}{\mathcal{E}}
\newcommand{\cF}{\mathcal{F}}
\newcommand{\cG}{\mathcal{G}} 
\newcommand{\cI}{\mathcal{I}} 
\newcommand{\cO}{\mathcal{O}} 
\newcommand{\cM}{\mathcal{M}} 
\newcommand{\cN}{\mathcal{N}} 
\newcommand{\cP}{\mathcal{P}} 
\newcommand{\cQ}{\mathcal{Q}} 
\newcommand{\cS}{\mathcal{S}} 
\newcommand{\cU}{\mathcal{U}} 
\newcommand{\cJ}{\mathcal{J}}
\newcommand{\cX}{\mathcal{X}}
\newcommand{\cT}{\mathcal{T}}
\newcommand{\cV}{\mathcal{V}}
\newcommand{\cW}{\mathcal{W}}
\newcommand{\cB}{\mathcal{B}}
\newcommand{\cR}{\mathcal{R}}
\newcommand{\cH}{\mathcal{H}}
\newcommand{\cZ}{\mathcal{Z}}
\newcommand{\D}{\bar{B}}

\newcommand{\Ker}{\textrm{Ker}\,}
\renewcommand{\coker}{\textrm{coker}\,}

\newcommand{\ext}{\mathrm{ext}} 
\newcommand{\x}{\times}

\newcommand{\mM}{\mathscr{M}} 

\newcommand{\CC}{\mathbb{C}} 
\newcommand{\QQ}{\mathbb{Q}} 
\newcommand{\FF}{\mathbb{F}} 
\newcommand{\PP}{\mathbb{P}} 
\newcommand{\HH}{\mathbb{H}} 
\newcommand{\RR}{\mathbb{R}} 
\newcommand{\ZZ}{\mathbb{Z}} 
\newcommand{\NN}{\mathbb{N}} 
\newcommand{\DD}{\mathbb{D}} 

\renewcommand{\lg}{\mathfrak{g}} 
\newcommand{\lh}{\mathfrak{h}} 
\newcommand{\lu}{\mathfrak{u}} 
\newcommand{\la}{\mathfrak{a}} 
\newcommand{\lb}{\mathfrak{b}} 
\newcommand{\lm}{\mathfrak{m}} 
\newcommand{\lgl}{\mathfrak{gl}} 
\newcommand{\lZ}{\mathfrak{Z}} 

\newcommand{\too}{\longrightarrow}
\newcommand{\imat}{\sqrt{-1}} 
\newcommand{\tinyclk}{{\scriptscriptstyle \Taschenuhr}} 
\newcommand\restr[2]{\left.#1\right|_{#2}}
\newcommand\rtorus[1]{{\mathbb{T}^{#1}}}
\newcommand\actPartial{\overline{\partial}}
\newcommand\handle[2]{\mathcal{A}^{#1}_{#2}}

\newcommand\pastingArea[2]{S^{#1-1} \times \bar{B}^{#2-#1}}
\newcommand\pastingAreaPlus[2]{S^{#1} \times \bar{B}^{#2-#1-1}}
\newcommand\coordvector[2]{\left.\frac{\partial}{\partial {#1}}\right|_{#2}}

\newcommand\Sets{\textbf{Set}}
\newcommand\Cat{\textbf{Cat}}
\newcommand\Top{\textbf{Top}}
\newcommand\TopHLC{\textbf{Top}_{hlc}}
\newcommand\TopS{\textbf{Top}_\star}
\newcommand\Diff{\textbf{Diff}}
\newcommand\Diffc{\textbf{Diff}_c}

\newcommand\CBord[3]{\mathbf{Bd}_{{#1 #3}}^{#2}}
\newcommand\Bord[1]{\CBord{#1}{}{}}
\newcommand\EBord[2]{\CBord{#1}{#2}{}}
\newcommand\Bordo[1]{\CBord{#1}{or}{}}
\newcommand\Bordp[1]{\mathbf{Bdp}_{{#1}}}
\newcommand\Bordpar[2]{\mathbf{Bd}_{{#1}}(#2)}
\newcommand\Bordppar[2]{\mathbf{Bdp}_{{#1}}(#2)}

\newcommand\Tubo[1]{\mathbf{Tb}_{#1}^0}
\newcommand\Tub[1]{\mathbf{Tb}_{#1}}
\newcommand\ETub[2]{\mathbf{Tb}_{#1}^{#2}}
\newcommand\Tubp[1]{\mathbf{Tbp}_{#1}}
\newcommand\Tubpo[1]{\mathbf{Tbp}_{#1}^0}
\newcommand\Tubppar[2]{\mathbf{Tbp}_{#1}(#2)}

\newcommand\Embc{\textbf{Emb}_c}
\newcommand\EEmbc[1]{\textbf{Emb}_c^{#1}}
\newcommand\Embpc{\textbf{Embp}_c}
\newcommand\Embparc[1]{\textbf{Emb}_c(#1)}
\newcommand\Embpparc[1]{\textbf{Embp}_c(#1)}

\newcommand\cSp{\cS_p}
\newcommand\cSpar[1]{\cS_{#1}}

\newcommand\Diffpc{\textbf{Diff}_c}

\newcommand\PVar[1]{\mathbf{PVar}_{#1}}
\newcommand\PVarC{\mathbf{PVar}_{\CC}}
\newcommand\Sr[1]{\mathrm{S}{#1}}

\newcommand\Bordpo[1]{\mathbf{Bdp}_{{#1}}^{or}}
\newcommand\ClBordp[1]{\mathbf{l}\CBord{#1}{}{}}
\newcommand\CTub[3]{\mathbf{Tb}_{{#1 #3}}^{#2}}
\newcommand\CTubp[1]{\CTub{#1}{}{}}
\newcommand\CTubpp[1]{\mathbf{Tbp}_{#1}{}{}}
\newcommand\CTubppo[1]{\mathbf{Tbp}_{#1}^0}
\newcommand\CClose[3]{\mathbf{Cl}_{{#1 #3}}^{#2}}
\newcommand\CClosep[1]{\CClose{#1}{}{}}
\newcommand\Obj[1]{\mathrm{Obj}(#1)}
\newcommand\Mor[1]{\mathrm{Mor}(#1)}
\newcommand\Vect[1]{{#1}\textrm{-}\mathbf{Vect}}
\newcommand\Vecto[1]{{#1}\textrm{-}\mathbf{Vect}_0}
\newcommand\Mod[1]{{#1}\textrm{-}\mathbf{Mod}}
\newcommand\Modt[1]{{#1}\textrm{-}\mathbf{Mod}_t}
\newcommand\Rng{\mathbf{Ring}}
\newcommand\Grp{\mathbf{Grp}}
\newcommand\Grpd{\mathbf{Grpd}}
\newcommand\Grpdo{\mathbf{Grpd}_0}
\newcommand\HS[1]{\mathbf{HS}^{#1}}
\newcommand\MHS[1]{\mathbf{MHS}}
\newcommand\PHS[2]{\mathbf{HS}^{#1}_{#2}}
\newcommand\MHSq{\mathbf{HS}}
\newcommand\Sch{\mathbf{Sch}}
\newcommand\Sh[1]{\mathbf{Sh}\left(#1\right)}
\newcommand\QSh[1]{\mathbf{QSh}\left(#1\right)}
\newcommand\Var[1]{\mathbf{Var}_{#1}}
\newcommand\Varrel[1]{\mathbf{Var}/{#1}}
\newcommand\KVarrel[1]{\K{\left(\mathbf{Var}/{#1}\right)}}
\newcommand\KoVarrel[1]{\Ko{\left(\mathbf{Var}/{#1}\right)}}
\newcommand\CVar{\Var{\CC}}
\newcommand\PHM[2]{\cM_{#1}^p(#2)}
\newcommand\MHM[1]{\cM_{#1}}
\newcommand\HM[2]{\textrm{HM}^{#1}(#2)}
\newcommand\HMW[1]{\textrm{HMW}(#1)}
\newcommand\VMHS[1]{VMHS({#1})}
\newcommand\geoVMHS[1]{VMHS_g({#1})}
\newcommand\goodVMHS[1]{\mathrm{VMHS}_0({#1})}
\newcommand\Par[1]{\mathrm{Par}({#1})}
\newcommand\K[1]{\mathrm{K}#1}
\newcommand\Ko[1]{\mathrm{\tilde K}{#1}}
\newcommand\KM[1]{\mathrm{K}\MHM{#1}}
\newcommand\KMo[1]{\mathrm{K}{\MHM{#1}}_0}
\newcommand\Ab{\mathbf{Ab}}
\newcommand\CPP{\cP\cP}
\newcommand\Bim[1]{{#1}\textrm{-}\mathbf{Bim}}
\newcommand\Span[1]{\mathrm{Span}({#1})}
\newcommand\Spano[1]{\mathrm{Span}^{op}({#1})}
\newcommand\GL[1]{\mathrm{GL}_{#1}}
\newcommand\SL[1]{\mathrm{SL}_{#1}}
\newcommand\PGL[1]{\mathrm{PGL}_{#1}}
\newcommand\Rep[1]{\mathfrak{X}_{#1}}

\newcommand\Repred[1]{\mathfrak{X}_{#1}^{r}}
\newcommand\Repirred[1]{\mathfrak{X}_{#1}^{ir}}
\newcommand\Repdiag[1]{\mathfrak{X}_{#1}^{d}}
\newcommand\Reput[1]{\mathfrak{X}_{#1}^{ut}}

\newcommand\Qtm[1]{\mathcal{Q}_{#1}}
\newcommand\sQtm[1]{\mathcal{Q}_{#1}^0}
\newcommand\Fld[1]{\mathcal{F}_{#1}}
\newcommand\Id{\mathrm{Id}}

\newcommand\DelHod[1]{e\left(#1\right)}
\newcommand\RDelHod{e}
\newcommand\eVect{\mathcal{E}}
\newcommand\e[1]{\eVect\left(#1\right)}
\newcommand\intMor[2]{\int_{#1}\,#2}

\newcommand\Dom[1]{\mathcal{D}_{#1}}

\newcommand\Xf[1]{{X}_{#1}}					
\newcommand\Xs[1]{\mathfrak{X}_{#1}}							

\newcommand\Xft[2]{\overline{{X}}_{#1, #2}} 
\newcommand\Xst[2]{\overline{\mathfrak{X}}_{#1, #2}}			

\newcommand\Xfp[2]{X_{#1, #2}}			
\newcommand\Xsp[2]{\mathfrak{X}_{#1, #2}}						

\newcommand\Xfd[2]{X_{#1; #2}}			
\newcommand\Xsd[2]{\mathfrak{X}_{#1; #2}}						

\newcommand\Xfm[3]{\mathcal{X}_{#1, #2; #3}}		
\newcommand\Xsm[3]{X_{#1, #2; #3}}					

\newcommand\XD[1]{#1^{D}}
\newcommand\XDh[1]{#1^{cr}}
\newcommand\XU[1]{#1^{UT}}
\newcommand\XP[1]{#1^{U}}
\newcommand\XPh[1]{#1^{\upsilon}}
\newcommand\XI[1]{#1^{\iota}}
\newcommand\XTilde[1]{#1^{ncr}}
\newcommand\Xred[1]{#1^{r}}
\newcommand\Xirred[1]{#1^{ir}}

\newcommand\Char[1]{\cR_{#1}}
\newcommand\Chars[1]{\cR_{#1}}
\newcommand\CharW[1]{\mathscr{R}_{#1}}

\newcommand\Ch[1]{\textrm{Ch}\,{#1}}
\newcommand\Chp[1]{\textrm{Ch}^+{#1}}
\newcommand\Chm[1]{\textrm{Ch}^-{#1}}
\newcommand\Chb[1]{\textrm{Ch}^b{#1}}
\newcommand\Der[1]{\textrm{D}{#1}}
\newcommand\Dp[1]{\textrm{D}^+{#1}}
\newcommand\Dm[1]{\textrm{D}^-{#1}}
\newcommand\Db[1]{\textrm{D}^b{#1}}

\newcommand\Gs{\cG}
\newcommand\Gq{\cG_q}
\newcommand\Gg{\cG_c}
\newcommand\Zs[1]{Z_{#1}}
\newcommand\Zg[1]{Z^{gm}_{#1}}
\newcommand\cZg[1]{\cZ^{gm}_{#1}}

\newcommand\RM[2]{R\left(\left.#1\right|#2\right)}
\newcommand\RMc[3]{R_{#1}\left(\left.#2\right|#3\right)}
\newcommand\set[1]{\left\{#1\right\}}
\newcommand{\Stab}{\textrm{Stab}\,} 
\renewcommand{\tr}{\textrm{tr}\,}             
\newcommand{\htr}{\textrm{tr}_0\,}             
\newcommand\EuChS{E}             
\newcommand\EuCh[1]{E\left(#1\right)}             

\newcommand{\Ann}{\textrm{Ann}\,}
\newcommand{\Rad}{\textrm{Rad}\,}  
\newcommand\supp[1]{\mathrm{supp}{(#1)}}
\newcommand\coh[1]{\left[H_c^\bullet\hspace{-0.05cm}\left(#1\right)\right]}
\newcommand\Bt{\Theta}
\newcommand\Be{B_e}
\newcommand\re{\textrm{Re}\,}
\newcommand\imag{\textrm{Im}\,}
\newcommand\Kahc{\textbf{K\"ah}_c}
\newcommand{\Ss}[1]{\cD_{#1}} 

\newcommand\Ccs[1]{C_{cs}(#1)}
\newcommand\can{\textrm{can}}
\newcommand\var{\textrm{var}}
\newcommand\Perv[1]{\textrm{Perv}(#1)}
\newcommand\DR[1]{\textrm{DR}{#1}}
\renewcommand\Gr[2]{\textrm{Gr}_{#1}^{#2}\,}
\newcommand\ChV{\textrm{Ch}\,}
\newcommand\VerD{^{\textrm{Ve}}\DD}
\newcommand\DHol[1]{\textrm{D}^b_{\textrm{hol}}(\cD_{#1})}
\newcommand\RegHol[1]{\Mod{\cD_{#1}}_{\textrm{rh}}}
\newcommand\Drh[1]{\textrm{D}^b_{\textrm{rh}}(\cD_{#1})}
\newcommand\Dcs[2]{\textrm{D}^b_{\textrm{cs}}({#1}; {#2})}
\newcommand{\rat}{\mathrm{rat}}
\newcommand{\dmod}{\mathrm{Dmod}}

\hyphenation{mul-ti-pli-ci-ty}

\hyphenation{mo-du-li}

\begin{abstract}
In this paper, we compute the virtual classes in the Grothendieck ring of algebraic varieties of $\SL{2}(\CC)$-character varieties over compact orientable surfaces with parabolic points of semi-simple type. When the parabolic punctures are chosen to be semi-simple non-generic, we show that a new interaction phenomenon appears generating a recursive pattern.
\end{abstract}
\null
\vspace{-1.1cm}
\maketitle

\vspace{-0.8cm}


\section{Introduction}
\let\thefootnote\relax\footnotetext{\noindent \emph{2010 Mathematics Subject Classification}. Primary:
 14C30. 
 Secondary:
 57R56, 
 14L24, 
 14D21. 

\emph{Key words and phrases}: TQFT, character varieties, Geometric Invariant Theory.}

Let $X$ be topological space with finitely generated fundamental group, $\pi_1(X)$, and let $G$ be a reductive algebraic group. The set of representations of $\pi_1(X)$ into $G$, $\rho: \pi_1(X) \to G$, has naturally a structure of algebraic variety, the so-called $G$-representation variety of $X$ and denoted $\Rep{G}(X)$. Moreover, the group $G$ itself acts on $\Rep{G}(X)$ by conjugation so we can consider the associated Geometric Invariant Theory (GIT) quotient, $\Char{G}(X) = \Rep{G}(X) \sslash G$, called the $G$-character variety of $X$. Since two representations are equivalent if and only if they are conjugated, the character variety $\Char{G}(X)$ is the moduli space of representations of $\pi_1(X)$ into $G$ \cite{Nakamoto}.

The topology and geometry of character varieties is an active research area. One of the the reasons of this interest is the celebrated non-abelian Hodge correspondence. It states that, if $X = \Sigma$ is a closed orientable surface and $G = \SL{n}(\CC)$, then the character variety of $\Sigma$ is diffeomorphic to the moduli space of rank $n$ vector bundles on $\Sigma$ with fixed determinant and equipped with a flat connection \cite{SimpsonI, SimpsonII}, and to the moduli space of rank $n$ and degree $0$ Higgs bundles on $\Sigma$ with fixed determinant \cite{Corlette:1988,Simpson:1992}. Despite that the three moduli spaces are naturally complex algebraic varieties, the correspondences are not holomorphic. This endows $\Char{\SL{n}(\CC)}(\Sigma)$ with three different complex structures that give rise to the first non-trivial example of a hyperk\"ahler manifold \cite{Hitchin:1992}.

For this reason, a thorough analysis of the natural algebraic structure on $\Char{G}(\Sigma)$ is needed. However, even in the simplest cases, the problem is very hard. The first approach was accomplished by Hausel and Rodr\'iguez-Villegas in \cite{Hausel-Rodriguez-Villegas:2008}. There, they introduced an arithmetic method that computes the $E$-polynomial of the $\GL{n}(\CC)$-character variety by counting its number of points in finite fields, in the spirit of the Weil conjectures. In \cite{Mereb}, Mereb extended the results to the case $G = \SL{n}(\CC)$. This strategy has also been exploited in \cite{Schiffmann:2016}, \cite{Mozgovoy-Schiffmann} and \cite{Mellit} at the side of the moduli space of Higgs bundles.

Despite of the power of the arithmetic method, it is a pure combinatorial approach that barely gives information about the underlying geometric structure of the character variety. In this way, Logares, Mu\~noz and Newstead in \cite{LMN} initiated a more geometric approach to the computation of $E$-polynomials of character varieties. The key idea of this paper is to chop $\Rep{G}(\Sigma)$ into simpler pieces for which the $E$-polynomial can be easily computed. Then, using the additivity of the $E$-polynomial, the polynomial of the whole space can be obtained by summing up all the contributions. Finally, they understood the identifications that appear in the GIT quotient to get the $E$-polynomial of $\Char{G}(\Sigma)$.

Using this method, for $G=\SL{2}(\CC)$, in \cite{LMN} it is explicitly computed the $E$-polynomial of the character variety over a surface of genus $1$ and $2$, in \cite{MM:2016} for genus $3$ and in \cite{MM} for arbitrary genus. Moreover, in \cite{Baraglia-Hekmati:2016}, using a mix between the arithmetic and the geometric method, it was computed the $E$-polynomial of character varieties over orientable surfaces for $G = \SL{2}(\CC), \SL{3}(\CC)$ and over non-orientable surfaces for $G= \SL{2}(\CC)$.

An even harder challenge appears when we consider a parabolic structure on $\Sigma$. Roughly speaking, it is given by a set of tuples $Q = \left\{(p_1, \lambda_1), \ldots, (p_s, \lambda_s)\right\}$, where $p_1, \ldots, p_s \in \Sigma$ is a collection of different marked points, called the punctures, and $\lambda_1, \ldots, \lambda_s \subseteq G$ is a collection of conjugacy classes of elements of $G$, called the holonomies. In this case a $Q$-parabolic representation is a representation $\rho: \pi_1\left(\Sigma - \set{p_1, \ldots, p_s}\right) \to G$ such that, if $\gamma_i$ is the loop around the puncture $p_i$, then $\rho(\gamma_i) \in \lambda_i$. The set of $Q$-parabolic representations also form an algebraic variety $\Rep{G}(\Sigma, Q)$, called the parabolic $G$-representation variety of $\Sigma$. The corresponding GIT quotient $\Char{G}(\Sigma, Q) = \Rep{G}(\Sigma, Q) \sslash G$ is the moduli space of $Q$-parabolic representations, known as the parabolic $G$-character variety. The non-abelian Hodge correspondence extends naturally to the parabolic setting to give diffeomorphisms with the moduli space of logarithmic flat connections and with the moduli space of parabolic Higgs bundles \cite{Simpson:parabolic}.

However, very little is known about the algebraic structure of character varieties in the parabolic case. One of the most important advances was done in \cite{Mellit:2017} and \cite{Mellit:2019} for $G=\SL{n}(\CC)$, where the $E$-polynomial is computed for parabolic structures of generic semi-simple type i.e.\ the holonomy conjugacy classes $\lambda_1, \ldots, \lambda_s \subseteq \SL{n}(\CC)$ are orbits of semi-simple elements lying in some Zariski open set of $\SL{n}(\CC)^s$. For general parabolic structures, in \cite{LM} the case of at most two punctures in $\SL{2}(\CC)$ is considered over elliptic curves. However, the arithmetic method is limited to deal with generic punctures. On the other hand, the geometric method is based on very subtle stratifications of the representation varieties that is not clear how to generalize to arbitrary many punctures.

In order to overcome this problem, in \cite{GPLM:2017} a new method was introduced based on Topological Quantum Field Theories (TQFTs) in the context of the PhD Thesis project of the author \cite{Gonzalez-Prieto:Thesis}. This method exploits the recursive nature of character varieties that is widely presented in the literature \cite{Carlsson-Rodriguez-Villegas, Diaconescu:2017, Hausel-Letellier-Villegas:2013, Mozgovoy:2012}. The key idea of this method the following. Let $\MHS{}$ be the category of mixed Hodge structures and let $\K{\MHS{}}$ be its associated Grothendieck ring (a.k.a.\ $K$-theory ring). Let us also consider $\Bordppar{n}{\Lambda}$ the category of $n$-dimensional bordisms of pairs with parabolic data in a collection $\Lambda$ of conjugacy classes of $G$ (see Section \ref{sec:TQFTs} for a precise definition). Then, in \cite{GPLM:2017}, we constructed a lax monoidal functor $\Zs{G}: \Bordppar{n}{\Lambda} \to \Mod{\K{\MHS{}}}$ computing the virtual Hodge structure (i.e.\ the image in $\K{\MHS{}}$) of parabolic $G$-representation varieties with holonomies in $\Lambda$. Recall that this means that, if $W$ is a closed connected $n$-dimensional manifold, $\star \in W$ is a basepoint and $Q$ is a parabolic structure on $W$ then, seen as a bordism $(W, \star, Q): \emptyset \to \emptyset$, we have that $\Zs{G}(W, \star, Q): \K{\MHS{}} \to \K{\MHS{}}$ satisfies $\Zs{G}(W, \star, Q)(1) = \coh{\Rep{G}(W, Q)}$.

As an application, in \cite{GP:2018a} we used this method to compute the virtual Hodge structure of parabolic $\SL{2}(\CC)$-representation varieties over orientable surfaces of arbitrary genus and any number of punctures with Jordan-type holonomy. For this purpose, we showed that all the computations of the lax monoidal TQFT can be performed within a finitely generate $\K{\MHS{}}$-module, $\cW$, called the core submodule. This simplifies the calculations since, in that case, the TQFT can be described explicitly by computing the images of finitely many elements. Using the results of \cite{GP:2018b} about stratification of GIT quotients, we translated these results to give the virtual Hodge structures of $\SL{2}(\CC)$-character varieties over surfaces of arbitrary genus an any number of punctures of Jordan type.

Nonetheless, if we consider semi-simple holonomies in the parabolic structure, the situation becomes much more involved. The most important problem is that the semi-simple punctures trigger the appearance of new generators of the submodule $\cW$, so it is no longer invariant under the TQFT. Moreover, if the punctures are not generic, an interaction phenomenon arises between these new generators. The aim of this paper is to explore this case.

In Section \ref{sec:TQFTs}, we will sketch briefly the construction of \cite{GPLM:2017} of the lax monoidal TQFT, $\Zs{G}$. Moreover, we will show that, with a small modification, we can improve this TQFT to compute, not only virtual Hodge structures on representation varieties, but indeed the virtual class of the representation variety $[\Rep{G}(X, Q)] \in \Ko{\Var{\CC}}$. Here $\Ko{\Var{\CC}}$ is the localization of the usual Grothendieck ring of complex algebraic varieties, $\K{\Var{\CC}}$, by a certain multiplicative set. This extension is compatible with the description of the TQFT in \cite{GPLM:2017} and those computations can be translated directly to this new context.

Section \ref{sec:parabolic-sl2c-repr} is the core of this paper. There, we perform the computation of the TQFT for $G=\SL{2}(\CC)$ and punctures of semi-simple type. For this purpose, we explicitly identify the new generators induced by the semi-simple punctures. In Section \ref{sec:computation-TQFT}, we perform the computation of the TQFT for the tube with a single puncture with semi-simple holonomy and we express the result in terms of the new generators. In Section \ref{sec:image-Tt-other-tubes}, we extend the computations of Sections 6.3 and 6.4 of \cite{GP:2018a} to the new set of generators, completing the explicit description of the TQFT. Finally, in Section \ref{sec:interation-phenomenon} we address the interaction phenomenon, giving rise to a combinatorial formula that shows how to modify the generic virtual class to deal with the case of non-generic punctures. These interaction phenomena are at the bottom of the reason why the arithmetic method breaks down when considering non-generic punctures. Therefore, as a consequence of these computations, in Theorem \ref{thm:result-complete} we obtain the following result.

\vspace{-0.06cm}
\begin{theorem*}
Let $\Sigma_g$ be the closed orientable genus $g$ surface and fix traces $t_1, \ldots, t_s \in \CC-\set{\pm 2}$, maybe non-generic. Write them as $t_i = \lambda_i + \lambda_i^{-1}$ for some $\lambda_i \in \CC^*-\set{\pm 1}$. Let $\alpha_+$ (resp.\ $\alpha_-$) be one half of the number of tuples $(\epsilon_1, \ldots, \epsilon_s) \in \set{\pm 1}^s$ such that $\lambda_1^{\epsilon_1}\cdots \lambda_s^{\epsilon_s} = 1$ (resp.\ such that $\lambda_1^{\epsilon_1}\cdots \lambda_s^{\epsilon_s} = -1$). Let $Q$ be a parabolic structure with $r$ punctures with holonomy $[J_+]$ and $s>0$ punctures with holonomies $\Ss{t_1}, \ldots, \Ss{t_s}$. Denote $q = [\CC] \in \K{\Var{\CC}}$ and let $\Ko{\Var{\CC}}$ be the localization of $\K{\Var{\CC}}$ with respect to the multiplicative set generated by $q, q+1$ and $q-1$. The virtual class of $\Rep{\SL{2}(\CC)}(\Sigma_g, Q)$ in $\Ko{\Var{\CC}}$ is
\begin{itemize}
	\item If $r>0$, then
\begin{align*}
	\left[\Rep{\SL{2}(\CC)}(\Sigma_g, Q)\right] =\, & q^{2g + s-1}(q - 1)^{2g + r-1}(q+1)\left(2^{2g+s-1} - 2^s + (q + 1)^{2g + r+s-2} \right)\\
	&+ \overline{\cI}_r(t_1, \ldots, t_s),
\end{align*}
where the interaction term is given by
\begin{align*}
\overline{\cI}_r(t_1, \ldots, t_s) =\, & q^{2g + s-1}(q - 1)^{2g + r-1}(\alpha_+ + \alpha_-)\bigg( 2^{2g} + 2^{2g}q  - 2q -2 \\
& \left. + (q + 1)^{2g + r} + (q+1)\left(1 - 2^{2g-1} - \frac{1}{2}(q + 1)^{2g + r-1}\right)\right) \\
& + q^{2g + s-1}(q - 1)^{2g + r}(q+1)\alpha_+.
\end{align*}
	\item If $r=0$, then
\begin{align*}
	\left[\Rep{\SL{2}(\CC)}(\Sigma_g, Q)\right] =\, & q^{2g+s-1}(q - 1)^{2g-1}(q + 1)(2^{2g+ s - 1} - 2^s + (q + 1)^{2g+s-2} \\
	&+ q^{2-2g-s}(q + 1)^{2g+s-2}) + \overline{\cI}_0(t_1, \ldots, t_s),
\end{align*}
where the interaction term is given by
\begin{align*}
\overline{\cI}_0(t_1, \ldots, t_s) =\, & q^{s-1}(q - 1)^{2g-1}(q+1)(\alpha_+ + \alpha_-)\left(q(q + 1)^{2g-1} + q^{2g}(q + 1)^{2g-1}\right.
\\ 
&\left.- q^{2g}(q + 1)^{2g-1} - q(q + 1)^{2g-1}\right) +q^{2g+s-1}(q - 1)^{2g}(q + 1)\alpha_+.
\end{align*}
\end{itemize}
\end{theorem*}

Finally, in Section \ref{sec:character-varieties} we use the techniques of \cite{GP:2018b} to translate this result across the GIT quotient down to character varieties. For this reason, in Section \ref{sec:review-pseudoquotients} we review the theory of pseudo-quotients that allows stratifications of representation varieties. In Section \ref{sec:semi-simple-quotient} we apply this theory to count the identifications that take place in the GIT quotient of the parabolic representation variety. In this way, we finally obtain the following result (Theorem \ref{thm:result-complete-quotient}).

\begin{theorem*}
The virtual class of $\Char{\SL{2}(\CC)}(\Sigma_g, Q)$ in $\Ko{\Var{\CC}}$ is
\begin{itemize}
	\item If $r>0$, then
\begin{align*}
	\left[\Char{\SL{2}(\CC)}(\Sigma_g, Q)\right] =\, & q^{2g + s-2}(q - 1)^{2g + r-2}\left(2^{2g+s-1} - 2^s + (q + 1)^{2g + r+s-2} \right) + \frac{\overline{\cI}_r( t_1, \ldots, t_s)}{q^3-q}.
\end{align*}
	\item  If $r = 0$, then
\begin{align*}
	\left[\Char{\SL{2}(\CC)}(\Sigma_g, Q)\right] =\, & q^{2g + s - 2}(q - 1)^{2g - 2}\left(2^{2g + s - 1}- 2^s + (q + 1)^{2g + s - 2}\right) + (q^2 - 1)^{2g - 2}(q + 1)^{s} \\&+ 2\alpha_{+}(q - 1)^{2g-1}(2q - 2q^{2g + s - 2}-1) + \frac{\overline{\cI}_0( t_1, \ldots, t_s)}{q^3-q}.
\end{align*}
\end{itemize}
\end{theorem*}

Also, we will show how the remaining combinations of holonomies for the punctures can be reduced to one of these cases of the ones studied in \cite{GP:2018a, GP:2018b}.

This result finishes the study of virtual classes of parabolic $\SL{2}(\CC)$-character varieties. However, much remain to be done in this business. First, the next goal would be to extend these results to higher rank. The case $G=\SL{3}(\CC)$ would allow us to extend the results of \cite{Baraglia-Hekmati:2016} to the parabolic case and the case $G=\SL{4}(\CC)$ is completely unknown. We expect that these situations might be addressed with the techniques developed in this series of papers. Furthermore, it would be interesting to consider other families of groups as $\PGL{n}(\CC)$ or $\textrm{Sp}_n(\CC)$ to explore the similarities and differences in the corresponding TQFTs. This study will be relevant towards the understanding of the mirror symmetry conjectures for character varieties, since $\SL{n}(\CC)$ and $\PGL{n}(\CC)$ are Langlands dual groups.

The next step would be to extend the TQFT used here to deal with representation varieties over more general spaces. In particular, it would be interesting to consider the case of character varieties over singular and non-orientable surfaces, as well as over complements of knots. This is the objective of a upcoming paper.

A more ambitious goal would be to extend the TQFT across the non-abelian Hodge correspondence to compute also virtual classes of moduli spaces of flat connections and moduli spaces of Higgs bundles. This would allow us to capture not only a particular complex structure but the whole picture of the hyperk\"ahler structure. Finally, we expect that the TQFT constructed will be useful to shed some light into the mirror symmetry conjectures for character varieties posed in \cite{Hausel:2005} that predict some astonishing symmetries of $E$-polynomials of character varieties over a group $G$ and its Langlands dual group $^{L}G$.

\subsection*{Acknowledgements}

The author wants to thank Tam\'as Hausel, Anton Mellit, Peter Newstead and Thomas Wasserman for very useful conversations. I would also like to express my highest gratitude to my PhD advisors Marina Logares and Vicente Mu\~noz for their invaluable help, support and encouragement throughout the development of this paper.

The author acknowledges the hospitality of the Faculty of Mathematical Sciences at Universidad Complutense de Madrid, in which part of this work was completed. The author has been supported by a "\!la Caixa" Foundation scholarship LCF/BQ/ES15/10360013 and partially by a MINECO (Spain) Project MTM2015-63612-P.

\section{Topological Quantum Field Theory for representation varieties}
\label{sec:TQFTs}
In this section, we shall sketch briefly the construction of the TQFT described in \cite{GPLM:2017}. We will also include some modifications that will allow us to compute, not only virtual Hodge structures, but the whole virtual class in the Grothendieck ring of algebraic varieties.

We will follow the notation of \cite{GPLM:2017}. Let $\Lambda$ be a set, $n \geq 1$ and let $\Bordppar{n}{\Lambda}$ be the category of $n$-bordisms of pairs with parabolic data in $\Lambda$. Recall that an object of this category is a triple $(M, A, Q)$ with $M$ a closed $(n-1)$-dimensional manifold, $A \subseteq M$ a finite set meeting each connected component of $M$ and $Q = \left\{(S_1, \lambda_1), \ldots, (S_r, \lambda_r)\right\}$ a parabolic structure on $M$, i.e.\ a set of co-oriented disjoint submanifolds $S_1, \ldots, S_r \subseteq M$ of codimension $2$, with labels $\lambda_1, \ldots, \lambda_r \in \Lambda$.

A morphism $(W, A, Q): (M_1, A_1, Q_1) \to (M_2, A_2, Q_2)$ is given by a unoriented bordism $W$ between $M_1$ and $M_2$, a finite set of points $A \subseteq W$ meeting each connected component of $W$ and such that $M_1 \cap A=A_1$ and $M_2 \cap A=A_2$, and $Q$ a parabolic structure on $W$ such that the restrictions on $M_1$ and $M_2$ agree with $Q_1$ and $Q_2$, respectively. Composition of morphisms is given by gluing of bordisms along their common boundary and juxtaposition of basepoints and parabolic structures. The category $\Bordppar{n}{\Lambda}$ is naturally a monoidal category with monoidal product the usual disjoint union of manifolds.

Now, let $R$ be a commutative and unitary ring, and let $\Mod{R}$ be the category of $R$-modules and $R$-module homomorphisms. A \emph{lax monoidal Topological Quantum Field Theory}, shortened TQFT, is a lax monoidal functor
$$
	Z: \Bordppar{n}{\Lambda} \to \Mod{R}.
$$

\begin{rmk}
The lax monoidality condition means that, for any $(M_1, A_1, Q_1),(M_2, A_2, Q_2) \in \Bordppar{n}{\Lambda}$, there exists a $R$-module homomorphism
$$
	Z(M_1, A_1, Q_1) \otimes Z(M_2, A_2, Q_2) \to Z\left((M_1, A_1, Q_1) \sqcup (M_2, A_2, Q_2)\right).
$$
However, this morphism might be not an isomorphism, in contrast to what is mandatory for a genuine monoidal functor.
\end{rmk}

\begin{rmk}
We can endow $\Bordppar{n}{\Lambda}$ and $\Mod{R}$ with natural $2$-category structures. In this framework, a lax monoidal TQFT can be usually promoted to a $2$-functor. However, for computational purposes, we will not need this structure, so we will not explore it further in this paper. For more information, see \cite{GP:2018a}.
\end{rmk}

A closed $n$-dimensional manifold $W$, together with a finite set $A \subseteq W$ and a parabolic structure $Q$ on it, defines a morphism $(W,A,Q): \emptyset \to \emptyset$. Hence, under the TQFT, it gives rise to a $R$-linear map $Z(W,A,Q): R \to R$. This map is fully determined by the element $\chi(W,A,Q)=Z(W,A,Q)(1) \in R$ that we can think as an algebraic invariant of $(W,A,Q)$. In this sense, we say that $Z$ computes the invariant $\chi$.

\subsection{Standard TQFT}

In \cite{GPLM:2017}, it was constructed a lax monoidal TQFT that computes virtual Hodge structures of representation varieties. Such TQFT was constructed by means of a `pull-push construction' by splitting the TQFT into a `field theory' and a `quantisation', being the later constructed via a $\cC$-algebra (for a review of this method and these concepts, see \cite[Section 4]{GP:2018a}). In this section, we will slightly  extend this construction to compute the whole virtual class of the representation variety in the Grothendieck ring of algebraic varieties.

Let us fix a ground field $k$. Let $\Var{k}$ be the category of algebraic varieties over $k$ with regular morphisms between them. Moreover, given $X \in \Var{k}$, we will denote by $\Varrel{X}$ the relative category of algebraic varieties over $X$. Recall that the objects of this category are pairs $(Y,\pi)$ with $Y$ and algebraic variety and $\pi$ a morphism $\pi: Y \to X$. If the morphism $\pi$ is clear from the context, the object will be denoted just by $Y \in \Varrel{X}$. Given objects $(Y_1, \pi_1), (Y_2, \pi_2) \in \Varrel{X}$, a morphism between them is a regular morphism $f: Y_1 \to Y_2$ such that $\pi_1 = \pi_2 \circ f$. Finally, we will consider the associated Grothendieck rings (also known as the $K_0$-ring in $K$-theory), $\K{\Var{k}}$ and $\KVarrel{X}$. The image of an algebraic variety $X \in \Var{k}$ in the Grothendieck ring will be denoted $[X] \in \K{\Var{k}}$ and will be called the \emph{virtual class} of $X$.

Now, fix an algebraic group $G$ and let $\Lambda$ be a collection of subvarieties of $G$ that are invariant under conjugation (e.g.\ conjugacy classes of some elements). Given $n \geq 1$, we are going to construct a lax monoidal TQFT$, \Zs{G}: \Bordppar{n}{\Lambda} \to \Mod{\K{\Var{k}}}$,
such that, for all morphism $(W, A, Q): \emptyset \to \emptyset$, it gives $$
	\Zs{G}(W,A,Q)(1) = \left[\Rep{G}(W,A,Q)\right].
$$
This TQFT is called the \emph{standard TQFT}. In order to construct this functor, we are going to split it into two functors
$$
	\Bordppar{n}{\Lambda} \stackrel{\Fld{}}{\longrightarrow} \Span{\Var{k}} \stackrel{\Qtm{}}{\longrightarrow} \Mod{\K{\Var{k}}}.
$$
Here, $\Span{\Var{k}}$ is the category of spans of $\Var{k}$ (see \cite{Benabou} for the definition). The functor $\Fld{}$ is playing the role of a field theory and $\Qtm{}$ is playing the role of a quantisation (in the physical sense).

The field theory, $\Fld{}: \Bordppar{n}{\Lambda} \to \Span{\Var{k}}$, coincides with the one described in \cite{GPLM:2017}. On an object $(M, A, Q) \in \Bordppar{n}{\Lambda}$, it assigns $\Fld{}(M,A,Q) = \Rep{G}(M,A,Q)$, the associated parabolic $G$-representation variety of the fundamental groupoid $\Pi(M,A)$ (see Remark \ref{rmk:repr-variety}). Also, given a bordism $(W, A, Q): (M_1, A_1,Q_1) \to (M_2, A_2, Q_2)$, it assigns the span
$$
	\Rep{G}(M_1, A_1,Q_1) \stackrel{\,\,i_1}{\longleftarrow} \Rep{G}(W, A,Q) \stackrel{i_2}{\longrightarrow} \Rep{G}(M_2, A_2, Q_2).
$$
Here, $i_1$ and $i_2$ are the maps induced by the inclusions $M_1 \hookrightarrow W$ and $M_2 \hookrightarrow W$, respectively, at the level of representations. By the Seifert-van Kampen Theorem for fundamental groupoids, $\Fld{}$ is a well-defined functor.

\begin{rmk}\label{rmk:repr-variety}
Recall that, given a compact manifold $W$ and a parabolic structure $Q$, the set of $Q$-parabolic representations, $\rho: \pi_1(W) \to G$, has naturally the structure of an algebraic variety. Suppose that $Q=\left\{(S_1, \lambda_1), \ldots, (S_s, \lambda_s)\right\}$ is the parabolic structure. Roughly speaking, the algebraic structure comes from considering a finite set of generators $\gamma_1, \ldots, \gamma_r, \alpha_1, \ldots, \alpha_s$ of $\pi_1(W)$, with the loops $\alpha_i$ around the submanifold $S_i$ in the positive direction given by the orientation of the normal bundle of $S_i$. Then, we identify $\Rep{G}(W,Q)$ with the image of the map $\psi: \Rep{G}(W,Q) \to G^r \times \lambda_1 \times \ldots \times \lambda_s$, $\psi(\rho)=(\rho(\gamma_1), \ldots, \rho(\gamma_s), \rho(\alpha_1), \ldots, \rho(\alpha_r))$, that is an algebraic set.

In the case of having a set of basepoints $A \subseteq W$, we can consider the representation variety $\Rep{G}(W,A,Q)$ of representations of the fundamental groupoid $\rho: \Pi(W,A) \to G$ such that $\rho(\alpha_i) \in \lambda_i$ if $\alpha_i$ is a positive loop around $S_i$. In that case, by picking a distinguished element on each component of $\Pi(W,A)$, we have a natural identification $\Rep{G}(W,A,Q) = \Rep{G}(W_1,Q) \times \ldots \times \Rep{G}(W_m,Q) \times G^{|A|-m}$, where $W_1, \ldots, W_m$ are the connected components of $W$ meeting a basepoint. This also endows $\Rep{G}(W,A,Q)$ with the structure of an algebraic variety. This is the structure considered along this paper. 
\end{rmk}

The quantisation part $\Qtm{}: \Span{\Var{k}} \to \Mod{\K{\Var{k}}}$ is slightly different from the one in \cite{GPLM:2017}. We will construct it by means of a $\Var{k}$-algebra (see \cite[Section 4.1]{GP:2018a}) given by the following data.
\begin{itemize}
	\item On an object $X \in \Var{k}$, it assigns $\KVarrel{X}$, the Grothendieck ring of algebraic varieties relative to $X$. It is a ring with product the fibered product of varieties.
	\item In particular, the image of the singleton variety $\star \in \Var{k}$ (i.e.\ of the final object) is $\KVarrel{\star} = \K{\Var{k}}$. In this way, $\KVarrel{X}$ also has a natural $\K{\Var{k}}$-module structure by cartesian product.
	\item Given a regular morphism $f: X \to Y$, we define $f_!: \KVarrel{X} \to \KVarrel{Y}$ to be morphism given by $f_![Z, \pi] = [Z, f \circ \pi]$ for $[Z, \pi] \in \KVarrel{X}$. Analogously, we define $f^*: \KVarrel{Y} \to \KVarrel{X}$ by $f^*[Z, \pi] = [Z \times_Y X, \pi']$, where $\pi'$ is the map fitting in the pullback diagram
\[
   \xymatrix
   {	
	Z \times_Y X \ar[r]\ar[d]_{\pi'} & Z \ar[d]^{\pi} \\
	X \ar[r]_{f} & Y
   }
\]
Observe that $f_!$ is a $\K{\Var{k}}$-module homomorphism by its very definition, as it commutes with cartesian product. Moreover, $f^*$ is a ring homomorphism since $\left(Z_1 \times_Y Z_2\right) \times_Y X \cong \left(Z_1 \times_Y X\right) \times_X \left(Z_2 \times_Y X\right)$, for any $Z_1, Z_2 \in \Varrel{Y}$, and analogously for the corresponding regular morphisms.
\end{itemize}

By the usual base change property for algebraic varieties, this assignment has the Beck-Chevally property. Thus, it is a $\Var{k}$-algebra that we shall denote $\KVarrel{}$. Therefore, by \cite[Theorem 4.13]{GP:2018a}, it gives rise to a lax monoidal functor $\Qtm{} = \Qtm{\KVarrel{}}: \Span{\Var{k}} \to \Mod{\K{\Var{k}}}$. By construction, to a span of the form $S: \star \leftarrow X \rightarrow \star$, it assigns the homomorphism $\Qtm{}(S): \K{\Var{k}} \to \K{\Var{k}}$ such that $\Qtm{}(S)(1) = [X]$. Thus, putting together this quantisation functor with the field theory, we have obtained the following result.

\begin{thm}
Let $G$ be an algebraic group and $n \geq 1$. There exists a lax monoidal TQFT
$$
	\Zs{G}: \Bordppar{n}{\Lambda} \to \Mod{\K{\Var{k}}},
$$
computing virtual classes of parabolic $G$-representation varieties. 
\end{thm}

\subsection{Recovering Hodge monodromy representations}
Suppose that the ground field is $k=\CC$. In that case, we can also consider the $\Var{\CC}$-algebra of the $K$-theory of mixed Hodge modules, $\K{\MHM{}}$, as in \cite{GPLM:2017, GP:2018a}. In order to avoid misinterpretation with the maps of $\KVarrel{}$, the induced maps of $\K{\MHM{}}$ will be denoted $\cM f_!: \K{\MHM{X}} \to \K{\MHM{Y}}$ and $\cM f^*: \K{\MHM{Y}}\to \K{\MHM{X}}$. Note that this slightly differs from the notation of \cite{GPLM:2017, GP:2018a}.

Using $\K{\MHM{}}$ as the $\Var{\CC}$-algebra for the quantisation, we also obtain a lax monoidal TQFT that we will denote $\Zs{G}^{\cM}: \Bordppar{n}{\Lambda} \to \Mod{\K{\MHS{}}}$. Here $\MHS{} = \MHM{\star}$ is the category of rational mixed Hodge modules. This was the TQFT used in \cite{GP:2018a} for calculations.

This functor $\Zs{G}^{\cM}$ and the functor $\Zs{G}$ constructed here are strongly related.
To state the relation properly, suppose that we are working in a category $\cC$ with pullbacks and final object. Let $\cA=(A,B)$ and $\cA'=(A',B')$ be two $\cC$-algebras. In particular, this means that $A, A'$ are contravariant functors out of $\cC$ with values in rings, and $B, B'$ are covariant functors out of $\cC$ with values in modules over the rings $A(\star),A'(\star)$, respectively. By a natural transformation $\tau: \cA \Rightarrow \cA'$, we will refer to a collection of ring homomorphisms $\tau_c: \cA_{c} \to \cA'_{c}$, for $c \in \cC$, intertwining with the induced maps. This means that, for any morphism $f: c \to d$ of $\cC$, we have $\tau_c \circ A(f) = A'(f) \circ \tau_d$ and $\tau_d \circ B(f) = B'(f) \circ \tau_c$.

\begin{prop}\label{prop:relation-KMHM-KVar}
Let $X \in \Var{\CC}$. Consider the morphism $\tau_X: \KVarrel{X} \to \K{\MHM{X}}$ that, for $(Z, \pi) \in \Varrel{X}$, sends $\tau_X[Z, \pi] = \cM\pi_!\left(\underline{\QQ}_{Z}\right) \in \K{\MHM{X}}$, where $\underline{\QQ}_{Z} \in \K{\MHM{Z}}$ denotes de unit of the ring. Then $\tau$ defines a natural transformation of $\Var{\CC}$-algebras
$$
	\tau: \KVarrel{} \Rightarrow \K{\MHM{}}.
$$
\begin{proof}
The maps $\tau_X$ intertwine with the induced maps of the $\Var{\CC}$-algebras. In order to check it, let $f: X \to Y$ be a regular morphism. For the pushout, we directly have $\tau_Y \circ f_![Z,\pi] = \tau_Y[Z, f \circ \pi]=\cM(f \circ \pi)_! \underline{\QQ}_Z = \cM f_! \circ \cM\pi_! \left(\underline{\QQ}_Z\right) = \cM f_! \left(\tau_X[Z, \pi]\right)$. For the pullback map, consider the cartesian square
\[
   \xymatrix
   {	
	Z \times_Y X \ar[r]^{\;\;\;\;\;\;f'}\ar[d]_{\pi'} & Z \ar[d]^{\pi} \\
	X \ar[r]_{f} & Y
   }
\]

By the Beck-Chevalley property of $\K{\MHM{}}$, we have
\begin{align*}
	\tau_X \circ f^*[Z, \pi] &= \tau_X[Z \times_Y X, \pi'] = \cM\pi'_!\left( \underline{\QQ}_{Z \times_Y X}\right) = \cM\pi'_! \circ \cM f'^* \left(\underline{\QQ}_{Z}\right) \\
	&= \cM f^* \circ \cM\pi_! \left(\underline{\QQ}_{Z}\right) = \cM f^* \circ \tau_Y[Z, \pi].
\end{align*}
Observe that, in the third equality, we have used $\cM f'^* \left(\underline{\QQ}_{Z}\right) = \underline{\QQ}_{Z \times_Y X}$, since $\cM f'^*$ is a ring homomorphism.

Finally, let us show that the maps $\tau_X: \KVarrel{X} \to \K{\MHM{X}}$ are ring homomorphisms. Suppose that $\pi_1: Z_1 \to X$ and $\pi_2: Z_2 \to X$ are objects of $\Varrel{X}$. If $\Delta: X \to X \times X$ denotes the diagonal map, and $Z_1 \times Z_2$ denotes the usual cartesian product (i.e.\ the fibered product over $\star$), we have a cartesian square
\[
   \xymatrix
   {	
	Z_1 \times_X Z_2 \ar[r]^{\;\;\;\;\Delta'}\ar[d]_{} & Z_1 \times Z_2 \ar[d]^{\pi_1 \times \pi_2} \\
	X \ar[r]_{\Delta} & X \times X
   }
\]
By definition, we have $\cM\pi'_! \circ \cM\Delta'^*\left(\underline{\QQ}_{Z_1 \times Z_2}\right)=\cM\pi'_!\left(\underline{\QQ}_{Z_1 \times_X Z_2}\right) = \tau_{X}[Z_1 \times_X Z_2, \pi'] = \tau_{X}\left([Z_1, \pi_1] \cdot [Z_2, \pi_2]\right)$.

On the other hand, let us denote $p_1, p_2: X \times X \to X$ the projections onto the first and the second component, respectively. Since $ \underline{\QQ}_{Z_1} \boxtimes \underline{\QQ}_{Z_2} =\underline{\QQ}_{Z_1 \times Z_2}$, we have that
\begin{align*}
	\cM\Delta^* \circ& \cM(\pi_1 \times \pi_2)_!\left(\underline{\QQ}_{Z_1 \times Z_2}\right) = \cM\Delta^* \circ \cM(\pi_1 \times \pi_2)_!\left(\underline{\QQ}_{Z_1} \boxtimes \underline{\QQ}_{Z_2}\right)\\
	&= \cM\Delta^* \left(\cM (\pi_1)_!\underline{\QQ}_{Z_1} \boxtimes \cM (\pi_2)_!\underline{\QQ}_{Z_2}\right) = \cM\Delta^* \left(\tau_X[Z_1, \pi_1] \boxtimes \tau_X[Z_2, \pi_2]\right)\\
	&= \cM\Delta^* \left(\cM p_1^*\left(\tau_X[Z_1, \pi_1]\right) \otimes \cM p_2^*\left( \tau_X[Z_2, \pi_2]\right)\right) \\
	&= \cM (p_1 \circ \Delta)^*\tau_X[Z_1, \pi_1] \otimes \cM (p_2 \circ \Delta)^*\tau_X[Z_2, \pi_2] = \tau_X[Z_1, \pi_1] \otimes \tau_X[Z_2, \pi_2].
\end{align*}
Note that, in the last equality, we have used $p_1 \circ \Delta = p_2 \circ \Delta = 1_X$. Therefore, by the Beck-Chevalley property of $\K{\MHM{}}$, both elements agree. This proves that $\tau_X$ is a ring homomorphism.
\end{proof}
\end{prop}

\begin{rmk}
\begin{itemize}
\item The mixed Hodge module $\tau_X[Z,\pi] = \cM\pi_!\left(\underline{\QQ}_{Z}\right)$ first appeared in \cite{LMN} (see also \cite{MM}), where it was called the Hodge monodromy representation and was denoted by $R_\pi(Z)$, or $R(Z)$ if the map was clear from the context. In this notation, the intertwining property reads $\cM f^*R(Z)=R(f^*Z)$ and $\cM f_!R_{\pi}(Z) = R_{f \circ \pi}(Z)$. Moreover, the fact that $\tau_X$ is a ring homomorphism implies that $R(Z_1 \times_X Z_2)=R(Z_1) \otimes R(Z_2)$. 
\item The fact that $\tau_X$ is a ring homomorphism is quite surprising since, for a regular map $\pi: Z \to X$, the morphism $\cM\pi_!: \K{\MHM{Z}} \to \K{\MHM{X}}$ is not in general a ring homomorphism. However, it preserve the external product, that was exactly what we needed in order to complete the proof above.
\end{itemize}
\end{rmk}

In particular, the natural transformation $\tau$ gives us a ring homomorphism $\tau_\star: \K{\Var{\CC}} \to \K{\MHS{}}$. It induces a natural transformation $\Mod{\K{\MHS{}}} \Rightarrow \Mod{\K{\Var{\CC}}}$. Under this transformation, $\Zs{G}^{\cM}$ can be seen as taking values in $\Mod{\K{\Var{\CC}}}$. With these considerations, we have the following result.

\begin{cor}
There exists a natural transformation
$
	\Zs{G} \Rightarrow \Zs{G}^\cM,
$
given by $\tau_{\Rep{G}(M, A, Q)}: \KVarrel{\Rep{G}(M, A, Q)} \to \K{\MHM{\Rep{G}(M, A, Q)}}$, for $(M,A,Q) \in \Bordppar{n}{\Lambda}$.
\begin{proof}
By the construction of both functors, it is enough to build the natural transformation at the level of the respective quantisations, denoted $\Qtm{\KVarrel{}}, \Qtm{\K{\MHM{}}}: \Span{\Var{\CC}} \to \Mod{\K{\Var{\CC}}}$. For this purpose, we can just take $\tau_X: \Qtm{\KVarrel{}}(X) = \KVarrel{X} \to \Qtm{\K{\MHM{}}}(X)=\K{\MHM{X}}$, for an algebraic variety $X$. Unraveling the definitions, we get the claimed formula.
\end{proof}
\end{cor}

\subsection{Geometric and reduced TQFTs}

Despite that $\Zs{G}$ computes the virtual class of representations varieties, it is convenient to consider a modification of this TQFT that is easier to compute. It is produced by the reduction procedure described in \cite[Section 4.5]{GP:2018a}. In this Section, we will sketch the construction of this modification and we will describe the associated maps explicitly.

Given $(M, A, Q) \in \Bordppar{n}{\Lambda}$, there exists an action of $G$ on $\Rep{G}(M, A, Q)$ by conjugation. The orbit space of this action can be given the structure of a piecewise algebraic variety (see \cite[Section 5.2]{GP:2018a}), denoted $[\Rep{G}(M, A, Q) / G]$. We also have a piecewise algebraic quotient map $\pi: \Rep{G}(M, A, Q) \to [\Rep{G}(M, A, Q) / G]$.

With this piecewise quotient, we modify the field theory to assign, to any morphism $(W, A, Q): (M_1, A_1,Q_1) \to (M_2, A_2, Q_2)$ of $\Bordppar{n}{\Lambda}$, the span
$$
	[\Rep{G}(M_1, A_1,Q_1)/G] \stackrel{\pi \circ i_1}{\longleftarrow} \Rep{G}(W, A,Q) \stackrel{\pi \circ i_2}{\longrightarrow} [\Rep{G}(M_2, A_2, Q_2)/G].
$$
Using this new field theory and the quantisation induced by the $\Var{k}$-algebra $\KVarrel{}$, we obtain a new assignment
$
	\Zg{G}: \Bordppar{n}{\Lambda} \to \Mod{\K{\Var{k}}},
$
called the \emph{geometric TQFT}. 

Unfortunately, in this form $\Zg{G}$ is not a genuine functor since the new field theory does not satisfy the Seifert-van Kampen theorem. Nevertheless, this problem can be easily solved. Consider the endomorphism
$$
	\eta = \pi_! \circ \pi^*: \KVarrel{[\Rep{G}(M, A, Q)/G]} \to  \KVarrel{[\Rep{G}(M, A, Q)/G]}.
$$
If $\eta$ is invertible then, by \cite[Proposition 4.20]{GP:2018a}, $\cZg{G} = \Zg{G} \circ \eta^{-1}$ is a lax monoidal TQFT computing the same invariant than $\Zs{G}$, called the \emph{reduced TQFT}. For this reason, we can focus on the computation of the geometric TQFT, from which the reduced TQFT follows immediately.

\begin{rmk}\label{rmk:localization}
It may happen, and it will be actually the case in our computations of Section \ref{sec:parabolic-sl2c-repr}, that $\eta$ is not invertible directly. However, suppose that there exists a multiplicative set $S \subseteq \K{\Var{\CC}}$ such that $\eta: S^{-1}\left(\KVarrel{[\Rep{G}(M, A, Q)/G]}\right) \to S^{-1}\left(\KVarrel{[\Rep{G}(M, A, Q)/G]}\right)$ is actually invertible. In that situation, we can still get a reduced TQFT that, now, it will compute the virtual class of the representation variety in the localized ring $S^{-1}\left(\K{\Var{\CC}}\right)$.
\end{rmk}

\begin{figure}[h]
	\begin{center}
	\includegraphics[scale=0.47]{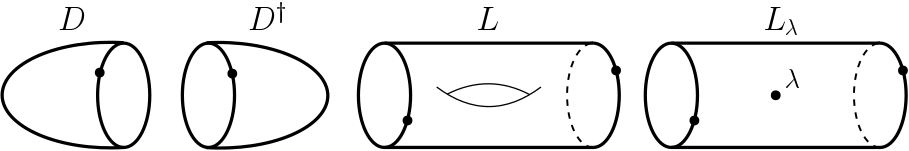}
	\end{center}
	\caption{Special bordisms of $\Bordppar{2}{\Lambda}$.}
	\label{fig:img-tubes}
\vspace{-0cm}
\end{figure}

In the case $n=2$, we can give explicitly the morphisms associated to the geometric TQFT. Observe that, in this case, the boundaries have no parabolic structures since the parabolic structure are codimension $2$ submanifolds. Now, consider the set of morphisms of $\Bordppar{2}{\Lambda}$ depicted in Figure \ref{fig:img-tubes}, for $\lambda \in \Lambda$.

The importance of these bordisms is the following. Let $\Sigma_g$ be the orientable closed surface of genus $g$ and consider any parabolic structure $Q$ on it with $s$ marked points. Then, if $A$ is a set of $g+s+1$ basepoints of $\Sigma_g$, as a morphism we can decompose $(\Sigma_g,A,Q) = D^\dag \circ L_{\lambda_s} \circ \ldots \circ L_{\lambda_1} \circ L^g \circ D$. Thus, in order to compute virtual images of representation varieties over closed surfaces, we can only focus on this special set of bordisms.

With respects to the boundaries, we only need to care about the circle with a single basepoint, $(S^1, \star)$. We have that the associated module is $\Zg{G}(S^1, \star) = \KVarrel{[\Rep{G}(S^1,\star)/G]} = \KVarrel{[G/G]}$ i.e.\ Grothendieck ring of varieties over the set of conjugacy classes of $G$. In order to get in touch with the notation of Section \ref{sec:parabolic-sl2c-repr}, let us denote this piecewise quotient map as $\htr: G \to [G/G]$. With respect to morphisms, as explained in \cite[Section 5.1]{GP:2018a}, the associated field theories of the discs are

\begin{minipage}[t]{0.5\textwidth}
$$
\begin{matrix}
	\Fld{}(D):& \star & \stackrel{}{\longleftarrow} & \star & \stackrel{i}{\longrightarrow} & [G/G], \\
	& &  & 1 & \mapsto & \set{\Id}
\end{matrix}
$$
\end{minipage}
\begin{minipage}[t]{0.5\textwidth}
$$
\begin{matrix}
	\Fld{}(D^\dag):& & [G/G] & \stackrel{i}{\longleftarrow} & \star & \stackrel{}{\longrightarrow} & \star. \\
	 & & \set{\Id} & \mapsfrom & 1 &  &
\end{matrix}
$$
\end{minipage}
Therefore, $\Zg{G}(D)=i_!$ and $\Zg{G}(D^\dag)=i^*$. On the other hand, for the tubes $L$ and $L_\lambda$, the associated field theories are
$$
\begin{matrix}
	\Fld{}(L): & [G/G] & \stackrel{p}{\longleftarrow} & G^4 & \stackrel{q}{\longrightarrow} & [G/G], \\
	& \htr (g) & \mapsfrom & (g, g_1, g_2, h) & \mapsto & \htr (g[g_1,g_2])
\end{matrix}
$$
$$
\begin{matrix}
	\Fld{}(L_\lambda):&[G/G] & \stackrel{r}{\longleftarrow} & G^2 \times \lambda & \stackrel{s}{\longrightarrow} & [G/G]. \\
	&\htr (g) & \mapsfrom & (g, h, \xi) & \mapsto & \htr (g\xi)
\end{matrix}
$$
Hence, $\Zg{G}(L)=q_!p^*$ and $\Zg{G}(L_\lambda)=s_!r^*$. Finally, the reduced TQFT is $\cZg{G} = \Zg{G} \circ \eta^{-1}$ with $\eta = \htr_! \circ \htr^*: \KVarrel{[G/G]} \to \KVarrel{[G/G]}$.

\section{Parabolic $\SL{2}(\CC)$-representation varieties}
\label{sec:parabolic-sl2c-repr}

From now on, we shall focus on the case of surfaces and $G = \SL{2}(\CC)$, for which we will compute the module homomorphisms of the reduced TQFT, as described in the previous section. This is analogous to the results of \cite{GP:2018a} but for a wider set of allowed holonomies for the punctures. For convenience, we will shorten $\Zg{\SL{2}(\CC)}$ by $\Zg{}$.

\subsection{Generalities on $\SL{2}(\CC)$}
In order to fix notation, recall that in $\SL{2}(\CC)$ there are five special types of elements, namely the matrices
$$
	\Id = \begin{pmatrix}
	1 & 0\\
	0 & 1\\
\end{pmatrix}, \hspace{0.2cm}
	-\Id = \begin{pmatrix}
	-1 & 0\\
	0 & -1\\
\end{pmatrix}, \hspace{0.2cm}
	J_+ = \begin{pmatrix}
	1 & 1\\
	0 & 1\\
\end{pmatrix}, \hspace{0.2cm}
	J_- = \begin{pmatrix}
	-1 & 1\\
	0 & -1\\
\end{pmatrix}, \hspace{0.2cm}
D_\lambda = \begin{pmatrix}
	\lambda & 0\\
	0 & \lambda^{-1}\\
\end{pmatrix},
$$
with $\lambda \in \CC^*-\left\{\pm 1\right\}$. Any element of $\SL{2}(\CC)$ is conjugated to one of these elements. Such a distinguished representant is unique up to the fact that $D_\lambda$ and $D_{\lambda^{-1}}$ are conjugated for all $\lambda \in \CC^*-\left\{\pm 1\right\}$. Hence, we have a stratification
\begin{equation*} \label{eqn:dec-sl2}
	\SL{2}(\CC) = \left\{\Id\right\} \sqcup \left\{-\Id\right\} \sqcup [J_+] \sqcup [J_-] \sqcup D,
\end{equation*}
where $
D = \bigcup_{\lambda} [D_\lambda] = \left\{A \in \SL{2}\,|\, \tr A \neq \pm 2\right\}$ and $[A]$ denotes the conjugacy class of $A \in \SL{2}(\CC)$. Given $t \in \CC-\set{\pm 2}$, we will also denote $\Ss{t} = \left\{A \in \SL{2}(\CC)\,|\, \tr A = t\right\}$.

The GIT quotient of the action of $\SL{2}(\CC)$ on itself is given by the trace map $\tr: \SL{2}(\CC) \to \CC$. By the stratification above, it only identifies the orbits $\set{\Id}$ and $[J_+]$ (both orbits of matrices of trace $2$) as well as the orbits $\set{-\Id}$ and $[J_-]$ (both orbits of matrices of trace $-2$). This implies that the piecewise quotient space is
\begin{equation*}
	[\SL{2}(\CC) / \SL{2}(\CC)] = \left\{[\Id]\right\} \sqcup \left\{[-\Id]\right\} \sqcup \set{[J_+]} \sqcup \set{[J_-]} \sqcup \Bt,
\end{equation*}
where $\Bt = \CC - \left\{\pm 2\right\}$ is the space of traces of matrices with two different simple eigenvalues. The orbitwise quotient map will be denoted $\htr: \SL{2}(\CC) \to [\SL{2}(\CC) / \SL{2}(\CC)]$. Under the action of $\SL{2}(\CC)$ on itself by conjugation we have that $\Stab(\pm \Id) = \SL{2}(\CC)$, $\Stab(J_\pm) = \CC$ and $\Stab(D_\lambda) = \CC^*$, for $\lambda \in \CC^*-\set{\pm 1}$. In particular, if we denote $q=[\CC] \in \K{\Var{\CC}}$, we have that $[\SL{2}(\CC)]=q^3-q$, $\left[[J_\pm]\right] = q^2-1$ and $[\Ss{t}] = q^2+q$ (see \cite[Section 6]{GP:2018a} for details).

In analogy with \cite[Section 6.1]{GP:2018a}, under the inclusion of these strata into $\SL{2}(\CC)$, we can consider the respective units $T_{2} \in \KVarrel{\set{\Id}}$, $T_{-2} \in \KVarrel{\set{-\Id}}$, $T_+ \in \KVarrel{[J_+]}$, $T_- \in \KVarrel{[J_-]}$ and $T_{\Bt} \in \KVarrel{\Bt}$ seen as elements of $\KVarrel{[\SL{2}(\CC)/\SL{2}(\CC)]}$.
On $\Bt$, we also consider the varieties
$$
	X_{2} = \left\{(w, t) \in \CC^* \times \Bt\,|\, w^2 = t - 2\right\}, \hspace{0.5cm} 	X_{-2} = \left\{(w, t) \in \CC^* \times \Bt\,|\, w^2 = t + 2\right\},
$$
Let $\pi$ be their projections onto the second component, $\pi(w,t) = t$. In this setting, the pairs $(X_{2}, \pi)$ and $(X_{-2}, \pi)$ define elements of $\KVarrel{\Bt}$ that we will denote $[X_2]$ and $[X_{-2}]$, respectively. We set $S_{2} = [X_{2}] - T_{\Bt}$ and $S_{-2} = [X_{-2}] - T_{\Bt}$. It can be proven that, if we consider $X_{2,-2} = \left\{(w, t) \in \CC^* \times \Bt\,|\, w^2 = t^2 - 4\right\}$, with projection over $\Bt$ given by $\pi(w,t) = t$, then $[X_{2,-2}] = T_{\Bt} + S_2 \times S_{-2}$, where $S_2 \times S_{-2}$ denotes the product in $\KVarrel{\Bt}$.
These elements satisfy some algebraic relations in $\KVarrel{\Bt}$ that we will need for the subsequent computations.

\begin{lem} \label{lem:relations-KVarTheta}
In the ring $\KVarrel{\Bt}$, the following relations hold.
$$
	S_2^2 = S_{-2}^2 = T_{\Bt}, \hspace{1cm} [X_{2,-2}] \times [X_{2}] = [X_{2,-2}] \times [X_{-2}] = [X_2 \times_{\Bt} X_{-2}].
$$
\begin{proof}
For the first equality, observe that it holds
\begin{align*}
	X_2 \times_{\Bt} X_2 &= \left\{(w_1, w_2, t)\,\left|\,w_1^2 = t-2, w_2^2 = t-2\right.\right\} \\
	&= \left\{(w_1, w_1, t)\,\left|\,w_1^2 = t-2\right.\right\} \sqcup \left\{(w_1, -w_1, t)\,\left|\,w_1^2 = t-2\right.\right\} \cong X_2 \sqcup X_2.
\end{align*}
Therefore, we have that $2T_{\Bt} + 2S_2 = 2[X_2] = [X_2] \times [X_2] = (T_{\Bt} + S_2)^2 = T_{\Bt} + 2S_{2} + S_2^2$ and, solving for $S_2^2$, we find $S_2^{2} = T_{\Bt}$. The computation for $S_{-2}^2$ is analogous.
The second equality can be obtained directly from this result, since we have
\begin{align*}
	[X_{2,-2}] \times [X_2] &= (T_{\Bt} + S_{2} \times S_{-2})(T_{\Bt} + S_2) = T_{\Bt} + S_2 + S_2^2 \times S_{-2} + S_2 \times S_{-2} \\
	&= T_{\Bt} + S_2 +  S_{-2} + S_2 \times S_{-2} = [X_2] \times [X_{-2}].
\end{align*} 
However, a more directed proof can be given. Using the explicit expression of $X_{2, -2} \times_{\Bt} X_2$, we have that
\begin{align*}
	X_{2,-2} \times_{\Bt} X_2 &= \left\{(w_1, w_2, t) \,\left|\, w_1^2 = (t-2)(t+2), w_2^2=t-2\right.\right\}\\
	& \cong \left\{(w_1', w_2, t) \,\left|\, w_1'^2 = t+2, w_2^2=t-2\right.\right\} = X_{2} \times_{\Bt} X_{-2}. 
\end{align*}
Observe that, in the second equality, the isomorphism is given by the map $(w_1, w_2, t) \mapsto (w_1/w_2, w_2, t)$.
\end{proof}
\end{lem}


The importance of these elements comes from the fact that, in \cite{GP:2018a}, it was proven that the submodule generated by $\cW = \langle T_{\pm 2}, T_\pm, T_{\Bt}, S_{\pm 2}, S_2 \times S_{-2} \rangle$ is invariant under $\Zg{}(L)$ and $\Zg{}(L_{[J_\pm]})$ and that $\Zg{}(D)$ and $\Zg{}(D^{\dag})$ are respectively the inclusion and projection onto $T_{2}$. This means that all the computations of the TQFT can be safely restricted to this submodule.

For the purposes of this paper, we also need to consider a kind of `skyscraper generators' over selected points of $\Bt$. For $t \in \Bt$, we shall denote by $T_{t}$ the image of the unit in $\KVarrel{\Ss{t}} \subseteq \KVarrel{\Bt}$. In this way, we will consider the $\K{\Var{\CC}}$-submodule generated by $\cW$ and all the skyscraper generators $T_t$ for $t \in \Bt$, that we will denote $\overline{\cW}$. Observe that $\overline{\cW}$ is not finitely generated, in contrast to $\cW$. 


\begin{rmk}
There is a small change of notation with respect to \cite[Section 6.1]{GP:2018a}. In that paper, the elements $T_2$ and $T_{-2}$ are denoted as $T_1$ and $T_{-1}$, respectively. We have decided to change notation in order to avoid the ambiguity of using $T_{\pm 1}$ to denote both the unit over $\pm \Id$ and over the matrices of trace $\pm 1$. With the new convention, $T_t$ always refers to the unit over $\Ss{t}$, the diagonalizable matrices of trace $t$.

Another minor change of notation is that $S_2 \times S_{-2}$ in this paper stands for what it was denoted as $S_2 \otimes S_{-2}$ in \cite{GP:2018a}. This change allows us to be consistent with the product notation in the quantisation ring, that is cartesian product of varieties in this paper but it was tensor product of mixed Hodge modules in \cite{GP:2018a}.
\end{rmk}

Let us consider the module endomorphism $\eta = \htr_! \circ \htr^*: \KVarrel{[\SL{2}(\CC)/\SL{2}(\CC)]} \to \KVarrel{[\SL{2}(\CC)/\SL{2}(\CC)]}$. The image under this map of $\cW$ was calculated in \cite{GP:2018a}. Taking into account the new skyscraper generators, we have the following result.

\begin{lem}\label{lem:eta-morph}
The submodule $\overline{\cW}$ is invariant for the morphism $\eta$. Indeed, the image of $T_{t}$ is given by $\eta(T_{t}) = (q^2 + q) T_{t}$ and, with respect to the standard set of generators, the matrix of $\eta$ is

$$\left(\begin{array}{c|ccccccccc}
& T_2 & T_{-2} & T_+ & T_- & T_{\Bt} & S_2 & S_{-2} & S_2 \times S_{-2} & T_t \\
\hline
T_2 & 1 & 0 & 0 & 0 & 0 & 0 & 0 & 0 & 0 
\\
T_{-2} & 0 & 1 & 0 & 0 & 0 & 0 & 0 & 0 & 0 
\\
T_{+} & 0 & 0 & q^{2} - 1 & 0 & 0 & 0 & 0 & 0 & 0 
\\
T_{-} & 0 & 0 & 0 & q^{2} - 1 & 0 & 0 & 0 & 0 &0 
\\
T_{\Bt} & 0 & 0 & 0 & 0 & q^{2} & 0 & 0 & q & 0 
\\
S_{2} & 0 & 0 & 0 & 0 & 0 & q^{2} & q & 0 & 0 
\\
S_{-2} & 0 & 0 & 0 & 0 & 0 & q & q^{2} & 0 & 0 
\\
S_2 \times S_{-2} & 0 & 0 & 0 & 0 & q & 0 & 0 & q^{2} & 0 
 \\
T_{t} & 0 & 0 & 0 & 0 & 0 & 0 & 0 & 0 & q^{2} +
q  \\
\end{array}\right)
$$
\begin{proof}
We just need to compute the image of the new skyscraper generators. Fix $t \in \Bt=\CC - \set{\pm 2}$. We have a commutative diagram whose square is cartesian
\[
\begin{displaystyle}
   \xymatrix
   {	
	& \Ss{t} \ar[r] \ar@{^{(}-{>}}[d] \ar[dl] & \left\{\Ss{t}\right\} \ar@{_{(}-{>}}[d] \\
	[\SL{2}(\CC)/\SL{2}(\CC)] & \SL{2}(\CC) \ar[r]^{\htr\;\;\;\;\;\;\;\;\;\;} \ar[l]_{\;\;\;\;\;\;\;\;\;\;\;\;\;\htr}& [\SL{2}(\CC)/\SL{2}(\CC)]
   }
\end{displaystyle}   
\]
Hence, we have that $\htr_! \circ \htr^* T_{t} = [\Ss{t}]\, T_t = (q^2 + q) T_t$.  Combining this computation with \cite[Proposition 6.3]{GP:2018a}, the result follows.
\end{proof}
\end{lem}

\begin{rmk}\label{rmk:loc-kvar}
The morphism $\eta: \KVarrel{[\SL{2}(\CC)/\SL{2}(\CC)]} \to \KVarrel{[\SL{2}(\CC)/\SL{2}(\CC)]}$ is not invertible, as the polynomials in $q$ appearing the the matrix above are not units of $\K{\Var{\CC}}$. However, $\eta$ is invertible as endomorphism of the localization of $\KVarrel{[\SL{2}(\CC)/\SL{2}(\CC)]}$ with respect to the multiplicative set $S$ generated by $q, q+1$ and $q-1$. From now on, we will denote the localizations $\KoVarrel{X} = S^{-1}\left(\KVarrel{X}\right)$ for a complex algebraic variety $X$ and, in particular, we set $\Ko{\Var{\CC}} = S^{-1}\left(\K{\Var{\CC}}\right)$. In this way, we get that $\eta: \KoVarrel{[\SL{2}(\CC)/\SL{2}(\CC)]} \to \KoVarrel{[\SL{2}(\CC)/\SL{2}(\CC)]}$ is invertible. As mentioned in Remark \ref{rmk:localization}, this implies we can still carry on our computations on the corresponding localized Grothendieck rings.
\end{rmk}

\subsection{The semi-simple puncture}

Fix $t_0 \in \Bt$ and write $t_0 = \lambda_0 + \lambda^{-1}_0$ for some $\lambda_0 \in \CC^* - \set{\pm 2}$. In this section, we shall compute the image of the tube $L_{\Ss{t_0}}: (S^1, \star) \to (S^1, \star)$ under the reduced TQFT, where $\Ss{t_0} = \left\{A \in \SL{2}(\CC)\,|\, \tr A = t_0\right\}$. Recall that the field theory for $\Zg{}$ on $L_{\Ss{t_0}}$ is the span
$$
\begin{matrix}
	[\SL{2}(\CC)/\SL{2}(\CC)] & \stackrel{{r}}{\longleftarrow} & \SL{2}(\CC)^2 \times \Ss{t_0} & \stackrel{{s}}{\longrightarrow} & [\SL{2}(\CC)/\SL{2}(\CC)]\\
	\htr(A) & \mapsfrom & (A, B, C) & \mapsto & \htr(BACB^{-1}) = \htr(AC)
\end{matrix}
$$
so $\Zg{}(L_{\Ss{t_0}}) = {s}_! \circ {r}^*: \KVarrel{[\SL{2}(\CC)/\SL{2}(\CC)]} \to \KVarrel{[\SL{2}(\CC)/\SL{2}(\CC)]}$.

Instead of the orbit space required by the field theory, it will be also convenient to consider the rough GIT quotient span
$$
\begin{matrix}
	\CC & \stackrel{}{\longleftarrow} & \SL{2}(\CC)^2 \times \Ss{t_0} & \stackrel{}{\longrightarrow} & \CC \\
	\tr(A) & \mapsfrom & (A, B, C) & \mapsto & \tr(AC)
\end{matrix}
$$
Taking into account that $\Ss{\lambda_0} = \SL{2}(\CC)/\CC^*$, we get that this span is naturally isomorphic to
$$
\begin{matrix}
	\CC & \stackrel{}{\longleftarrow} & \SL{2}(\CC)^2 \times \left(\SL{2}(\CC)/\CC^*\right) & \stackrel{}{\longrightarrow} & \CC \\
	\tr(A) & \mapsfrom & (A, B, P) & \mapsto & \tr(AD_{\lambda_0})
\end{matrix}
$$
Since the second and the third factors play no role in the previous span, we can focus on the simplified span
\begin{equation}
\label{eq:span-trace-row}
\begin{matrix}
	\CC & \stackrel{}{\longleftarrow} & \SL{2}(\CC) & \stackrel{}{\longrightarrow} & \CC \\
	\tr(A) & \mapsfrom & A & \mapsto & \tr(AD_{\lambda_0})
\end{matrix}
\end{equation}

\begin{lem}\label{lem:fibr-w}
Let $W$ be the algebraic variety 
$$
	W = \left\{(y,z,t,t') \in \CC^4\,\left|\,\,\frac{(t-t')^2}{t_0^2-4} - \frac{tt'}{t_0+2} +yz = -1\right.\right\}.
$$
Then, span (\ref{eq:span-trace-row}) is isomorphic to
$$
\begin{matrix}
	\CC & \stackrel{}{\longleftarrow} & W  & \stackrel{}{\longrightarrow} & \CC \\
	t & \mapsfrom & (y,z,t,t') & \mapsto & t'
\end{matrix}
$$
Moreover, let us decompose $W = W_{a} \sqcup W_{b}$ where $W_{a} = W \cap \left\{z \neq 0\right\}$ and $W_{b} = W \cap \left\{z = 0\right\}$. Then we have
\begin{itemize}
	\item $W_a$ is isomorphic to $\CC^* \times \CC \times \CC$ and the maps of the span correspond to the projections onto the second and third components, respectively.
	\item $W_b$ is isomorphic to $\CC \times \CC^*$. Under this isomorphism, the span is
$$
\begin{matrix}
	\CC & \stackrel{\alpha}{\longleftarrow} & \CC \times \CC^*& \stackrel{\beta}{\longrightarrow} & \CC \\
	\mu\lambda_0 + \mu^{-1}\lambda_0^{-1} & \mapsfrom & (y, \mu) & \mapsto & \mu + \mu^{-1}
\end{matrix}
$$
\end{itemize}
\begin{proof}
Let us take
$$
	W = \left\{(A, t, t') \in \SL{2}(\CC) \times \CC^2\,\left|\, \begin{matrix}t = \tr A\\ t' = \tr (AD_{\lambda_0})\end{matrix}\right.\right\}.
$$
A straightforward computation shows that $W$ is the set of matrices of the form
$$A=\begin{pmatrix}
	\frac{t}{2}\left(1-\frac{\lambda_0+\lambda_0^{-1}}{\lambda_0-\lambda_0^{-1}}\right) + \frac{t'}{\lambda_0-\lambda_0^{-1}} & y\\
	z & \frac{t}{2}\left(1+\frac{\lambda_0+\lambda_0^{-1}}{\lambda_0-\lambda_0^{-1}}\right)-\frac{t'}{\lambda_0-\lambda_0^{-1}} \\
\end{pmatrix}
$$	
satisfying that
$$
	\frac{\left(\lambda_0 t - \lambda_0 t'\right)^2}{(\lambda_0^2-1)^2} - \frac{\lambda_0 tt'}{(\lambda_0^2+1)^2} + yz = -1 \Leftrightarrow \frac{(t-t')^2}{t_0^2-4} - \frac{tt'}{t_0+2} +yz = -1.
$$
With this description, we have that the left-most span map corresponds to $ (y,z,t,t') \mapsto t$, and the right-most map to $(y,z,t,t') \mapsto t'$.
On $W_a$, we can solve this equation for $y$ and the span maps are trivial fibrations with fibers $\CC^* \times \CC$.

If $z = 0$, the later equation forces $t$ and $t'$ to satisfy
$$
H = \left\{(t,t') \in \CC^2 \left|\,\,\frac{(t-t')^2}{t_0^2-4} - \frac{tt'}{t_0+2} = -1\right.\right\}  = \left\{ t= \frac{1}{2} \, t' t_{0} \pm \frac{1}{2} \, \sqrt{{\left(t'^{2} - 4\right)}
{\left(t_{0}^{2} - 4\right)}}\right\}.
$$
Hence, $H$ is an affine hyperbola that can be parametrized by $\CC^* \to H$, $\mu \mapsto (\mu\lambda_0 + \mu^{-1}\lambda_0^{-1}, \mu + \mu^{-1})$. Under this isomorphism, the maps $r$ and $s$ become the claimed maps $\alpha$ and $\beta$, respectively.
\end{proof}
\end{lem}

\begin{prop} On $\Bt = \CC - \set{\pm 2}$, the endomorphism $\beta_! \circ \alpha^*: \KVarrel{\Bt} \to \KVarrel{\Bt}$ is given by
\begin{center}
\begin{tabular}{cc}
	$\beta_! \alpha^*(T_{\Bt}) = q\left(T_\Bt + S_2 \times S_{-2} - T_{t_0} - T_{-t_0}\right),$ & $\beta_! \alpha^*(S_2) = q\left(S_2 + S_{-2}- T_{t_0} - T_{-t_0}\right),$\\
	 $\beta_! \alpha^*(S_2 \times S_{-2}) = q \left(T_\Bt + S_2 \times S_{-2}- T_{t_0} - T_{-t_0}\right),$ & $\beta_! \alpha^*(S_{-2}) = q\left(S_2 + S_{-2}- T_{t_0} - T_{-t_0}\right).$\\
\end{tabular}
\end{center}
\begin{proof}
Observe that $\alpha^{-1}(\pm 2) = \CC \times \left\{\pm \lambda_0^{-1}\right\}$ and $\beta^{-1}(\pm 2) = \CC \times \left\{\pm 1\right\}$. Hence, the span for $W_b$ restricts to
$$
\begin{matrix}
	\CC-\set{\pm 2} & \stackrel{\alpha}{\longleftarrow} & \CC \times \left(\CC^* - \set{\pm 1, \pm \lambda_0^{-1}}\right)& \stackrel{\beta}{\longrightarrow} & \CC-\set{\pm 2} \\
	\mu\lambda_0 + \mu^{-1}\lambda_0^{-1} & \mapsfrom & (y, \mu) & \mapsto & \mu + \mu^{-1}
\end{matrix}
$$

Consider the isomorphism $\tau:  \CC^* - \set{\pm 1, \pm \lambda_0} \to  \CC^* - \set{\pm 1, \pm \lambda_0^{-1}}$ given by $\tau(\mu)=\lambda_0\mu$, the map $t(\mu) = \mu + \mu^{-1}$ and the projection $c: \CC \to \star$. We can write $\alpha = c \times (t \circ \tau)$ and $\beta = c \times t$ so we have that $\beta_!\alpha^* = c_!c^* \times t_!\tau^*t^* = q\,t_!\tau^*t^*$, with $q = [\CC] \in \K{\Var{\CC}}$.

In particular this implies that, for $T_{\Bt}$, it holds $\beta_!\alpha^*(T_{\Bt}) = q\,t_!t^*T_{\Bt}$ since $t^*$ and $\tau^*$ are ring homomorphisms. For $[X_2]$, we have that $t^*X_2 = \left(\CC^* - \set{\pm 1, \pm \lambda_0}\right) \times_{\Bt} X_2=\left\{(w, \mu)\,\left|\, w^2=\mu+\mu^{-1}-2\right.\right\}$, which is isomorphic to $(t \tau)^*X_2 = \left(\CC^* - \set{\pm 1, \pm \lambda_0^{-1}}\right) \times_{\Bt} X_2=\left\{(w, \mu)\,\left|\, w^2=\mu\lambda_0+\mu^{-1}\lambda_0^{-1}-2\right.\right\}$ as element of $\KVarrel{\left(\CC^* - \set{\pm 1, \pm \lambda_0^{-1}}\right)}$. Therefore, we have that $\beta_!\alpha^*[X_2] = q\,t_!t^*[X_2]$ and analogously for $X_2$ and $X_2 \times_{\Bt} X_{-2}$. Thus, $\tau^*$ does not modify the calculation and we can focus on the computation of $t_!t^*$. 

In this case, observe that we have an isomorphism $\psi: \CC^* - \set{\pm 1} \to X_{2,-2}$, given by $\psi(\mu) = (\lambda - \lambda^{-1}, \lambda + \lambda^{-1})$. Under this isomorphism, there is a commutative diagram 
	\[
\begin{displaystyle}
   \xymatrix
   {	
    & \CC^* - \set{\pm 1} \ar[rd]^{t} \ar[ld]_{t}& \\
	\CC - \set{\pm 2} & & \CC - \set{\pm 2} \\
	& X_{2,-2} \ar[ru]_{\pi} \ar[lu]^{\pi} \ar[uu]^{\psi} &
   }
\end{displaystyle}   
\]
where $\pi: X_{2,-2} \to \CC - \set{\pm 2}$ is $\pi(w,t)=t$. Since $\psi$ is an isomorphism, $t_!t^*=\pi_!\psi_!\psi^*\pi^* = \pi_!\pi^*$. In particular, we have that $\pi_!\pi^*(T_{\Bt}) = \pi_!\left(T_{X_{2,-2}}\right) = [X_{2,-2}] = T_{\Bt} + S_2 \times S_{-2}$. On the other hand, Lemma \ref{lem:relations-KVarTheta} implies that
$$
	\pi_!\pi^*[X_2] = [X_{2,-2} \times_{\Bt} X_{2}] = [X_2] \times [X_{-2}] = T_{\Bt} + S_{2} + S_{-2} + S_{2} \times S_{-2},
$$
and analogously for $[X_{-2}]$. For $[X_2 \times_{\Bt} X_{-2}]$, also using the relations of Lemma \ref{lem:relations-KVarTheta} we obtain
\begin{align*}
	\pi_!\pi^*[X_2 \times_{\Bt} X_{-2}] &= [X_{2,-2} \times_{\Bt} X_{2}] \times [X_{-2}] = [X_2] \times [X_{-2}] \times [X_{-2}] \\&= 2 [X_2] \times [X_{-2}] = 2T_{\Bt} + 2S_{2} + 2S_{-2} + 2S_{2} \times S_{-2}.
\end{align*}
Therefore, using that $[X_{\pm 2}] = T_{\Bt} + S_{\pm 2}$, and subtracting the missing fibers over $\pm t_0 = t(\pm \lambda_0^{-1})$, the result follows.
\end{proof}

\end{prop}

\section{Computations of the geometric TQFT}
\label{sec:computation-TQFT}

After the preliminary computations of Section \ref{sec:parabolic-sl2c-repr}, we are ready to accomplish the computation of the geometric TQFT. Along this section, we fix $t_0 \in \Bt$ and we write it as $t_0 = \lambda_0 + \lambda_0^{-1}$ for a fixed $\lambda_0 \in \CC^* - \set{\pm 1}$. We shall focus on the computation of the $\K{\Var{\CC}}$-module morphism 
$$
\Zg{}(L_{\Ss{t_0}}): \KVarrel{[\SL{2}(\CC) / \SL{2}(\CC)]} \to \KVarrel{[\SL{2}(\CC) / \SL{2}(\CC)]}.
$$

For short, if $X, Y$ are subvarieties of $[\SL{2}(\CC) / \SL{2}(\CC)]$, we shall denote
\begin{align*}
	Z_{X, Y} = {r}^{-1}(X) \cap s^{-1}(Y) &= \left\{(A, B, C) \in \htr^{-1}(X) \times \Ss{t_0} \,\,|\,\, \htr(AC) \in Y\right\} \\
	&\cong \left\{A \in \htr^{-1}(X) \,\,|\,\, \htr(AD_{\lambda_0}) \in Y\right\} \times \Ss{t_0} \times \SL{2}(\CC).
\end{align*}
We will also write $Z_{X,Y}^0 = \left\{A \in \htr^{-1}(X) \,\,|\,\, \htr(AD_{\lambda_0}) \in Y\right\}$ so $[Z_{X,Y}] = (q^2+q)(q^3-q)[Z_{X,Y}^0]$, where $q = [\CC] \in \K{\CVar}$.

\subsection{Image of $T_{2}$} We have that $r^{-1}(\Id) = Z_{\Id, \Ss{t_0}} = \Ss{t_0} \times \SL{2}(\CC)$ and $s$ is just the projection onto a point $s: \Ss{t_0} \times \SL{2}(\CC) \to \left\{\Ss{t_0}\right\}$. Hence, $\Zg{}(L_{\Ss{t_0}})(T_{2}) = (q^2+q)(q^3-q) T_{t_0}$.

\subsection{Image of $T_{-2}$} Analogously to the previous case, we have that $r^{-1}(-\Id) = \Ss{-t_0} \times \SL{2}(\CC)$ and $s$ is the projection onto the point $\set{\Ss{-t_0}}$ so $\Zg{}(L_{\Ss{t_0}})(T_{-2}) = (q^2+q)(q^3-q) T_{-t_0}$.

\begin{rmk}\label{rmk:trace-equal-0}
In general $D_{-\lambda_0}$ is not conjugated to $D_{\lambda_0}$ except in the case $\lambda_0 = \pm i$. In that case, $t_0 = 0$ so $\Zg{}(L_{\Ss{0}})(T_{2}) = \Zg{}(L_{\Ss{0}})(T_{-2}) = (q^2+q)(q^3-q) T_{0}$.
\end{rmk}

\subsection{Image of $T_+$} In this case, $r^{-1}([J_+])$ has a non-trivial decomposition. Analyzing the strata separately we obtain the following.
\begin{itemize}
	\item Over $\pm \Id$, we have that $Z_{[J_+], \set{\pm \Id}}^0 = \left\{A \sim J_+\,|\, AD_{\lambda_0} = \pm \Id\right\} = \emptyset$ so these strata add no contribution.
	\item Over $[J_\pm]$, a direct computation shows that
\begin{align*}
	Z_{[J_+], [J_\pm]}^0 &= \left\{A \sim J_+\,|\,AD_{\lambda_0} \sim J_\pm \right\} = \left\{\left.\begin{pmatrix}
	1-x & y\\
	z & 1+x\\
\end{pmatrix}\,\right| \begin{matrix}
	yz = -x^2\\
	\lambda_0 + \lambda_0^{-1} - x(\lambda_0 - \lambda_0^{-1}) = \pm 2\\
\end{matrix}
\right\} \\
	&= \left\{yz = -\left(\frac{\lambda_0 + \lambda_0^{-1}\mp 2}{\lambda_0 - \lambda_0^{-1}}\right)^2 \neq 0\right\} \cong \CC^*.
\end{align*}
This agrees with the calculation of Lemma \ref{lem:fibr-w}.
Therefore, $\Zg{}(L_{\Ss{t_0}})(T_{+})|_{[J_\pm]} = (q-1)(q^2+q)(q^3-q) T_{\pm}$.
	\item Over $\Bt$, by Lemma \ref{lem:fibr-w} we have that
$$
	Z_{[J_+], \Bt}^0 = \left\{\begin{matrix}\frac{\left(2\lambda_0 - \lambda_0 t'\right)^2}{(\lambda_0^2-1)^2} - \frac{2\lambda_0 t'}{(\lambda_0^2+1)^2} + yz = -1\\ (y,z) \neq (0,0),\,\,\,\,\,\,t' \neq \pm 2\end{matrix}\right\},
$$
and $s$ is given by $(y,z, t') \mapsto t'$. Where $z=0$, the previous equation forces $t' = t_0$. We have to remove the point $y=0$ from the fiber $\CC$, so this stratum contributes as $(q-1)T_{t_0}$. Where $z \neq 0$, we can solve for $y$ and $t' \in \Bt$ is free, so this stratum contributes as $(q-1)T_\Bt$. Hence, we have that $\Zg{}(L_{\Ss{t_0}})(T_{+})|_{\Bt} = (q-1)(q^2+q)(q^3-q)\left(T_{\Bt} + T_{t_0}\right)$.
\end{itemize}
Adding all the contributions of the strata, we finally have that
$$
\Zg{}(L_{\Ss{t_0}})(T_{+}) = (q-1)(q^2+q)(q^3-q)\left(T_+ + T_- + T_{\Bt} + T_{t_0}\right).
$$

\subsection{Image of $T_-$} The calculations for this stratum are completely analogous to the ones for $T_+$ so we find that
$$
\Zg{}(L_{\Ss{t_0}})(T_{-}) = (q-1)(q^2+q)(q^3-q)\left(T_+ + T_- + T_{\Bt} + T_{-t_0}\right).
$$
\subsection{Image of $T_{\Bt}$}\label{sec:image-T-Bt} As in the previous case, we stratify $r^{-1}(\Bt)$ and analyze each piece separately.
\begin{itemize}
	\item Over $\pm \Id$, we have $Z_{\Bt, \set{\pm \Id}} = \Ss{\pm t_0} \times \SL{2}(\CC)$ and $s$ is the projection onto a point. So we obtain that $\Zg{}(L_{\Ss{t_0}})(T_{\Bt})|_{\pm \Id} = (q^2+q)(q^3-q)T_{\pm 2}$.
	\item Over $[J_\pm]$, by Lemma \ref{lem:fibr-w}, we have that
\begin{align*}
	Z_{\Bt, [J_\pm]}^0 &= \left(W \cap \set{t'=2}\right) - \set{\pm D_{\lambda_0^{-1}}} \\&= \left\{(y,z,t)\,\left|\,\,\frac{\left(\lambda_0 t \mp 2\lambda_0\right)^2}{(\lambda_0^2-1)^2} \mp \frac{2\lambda_0 t}{(\lambda_0^2+1)^2} + yz = -1\right.\right\} - \set{(0,0, \pm t_0)}.
\end{align*}
Observe that we have to remove the point $A = \pm D_{\lambda_0^{-1}}$, for which $A D_{\lambda_0} = \pm\Id$ instead of $J_\pm$. This corresponds to $(y,z,t) = (0,0, \pm t_0)$. Where $z\neq 0$, we know by Lemma \ref{lem:fibr-w} that the fiber is $\CC^* \times \Bt$ that contributes with $(q-1)(q-2)$. Where $z = 0$, we have $t = \pm t_0$ so the fiber is $\CC^*$ that contributes with $q-1$. Taking into account the factor $(q^2+q)(q^3-q)$ we obtain that $\Zg{}(L_{\Ss{t_0}})(T_{\Bt})|_{[J_\pm]} = (q-1)^2(q^2+q)(q^3-q)T_{\pm}$.
	\item Over $\Bt$, analogously we have that
$$
	Z_{\Bt, \Bt}^0 = \left\{\frac{\left(\lambda_0 t - \lambda_0 t'\right)^2}{(\lambda_0^2-1)^2} - \frac{\lambda_0 tt'}{(\lambda_0^2+1)^2} + yz = -1\right\}.
$$
As proved in Lemma \ref{lem:fibr-w}, with respect to the projection onto $t'$, we have that for $z \neq 0$ it is a trivial fibration of fiber $\CC^* \times \Bt$ so this part contributes with $(q-1)(q-2)T_{\Bt}$. On the other hand, the stratum $z=0$ contributes with $q\left(T_{\Bt} + S_2 \times S_{-2} - T_{t_0} - T_{-t_0}\right)$.

\end{itemize}

Adding up the contributions, we find that
\begin{align*}
	\Zg{}(L_{\Ss{t_0}})(T_{\Bt})= \,& (q^2+q)(q^3-q)\left(T_2 + T_{-2} + (q-1)^2T_+ + (q-1)^2T_-\right.\\
	&+ \left.(q^2-2q+2) T_{\Bt} + q S_2 \times S_{-2} - q T_{t_0} - q T_{-t_0}\right).
\end{align*}

\subsection{Image of $S_{2}$} The computation for this generator is very similar the previous one but taking into account the double cover $X_2 \to \Bt$ that defines $S_2$. We shall focus on the computation of $s_!r^*[X_2]$ and, using that $[X_2] = T_{\Bt} + S_2$ the result follows by calculating $\Zg{}(L_{\Ss{t_0}})(S_2)=s_!r^*[X_2]-\Zg{}(L_{\Ss{t_0}})(T_{\Bt})$.

According to the stratification of $r^{-1}(\Bt)$ we have
\begin{itemize}
	\item Over $\pm \Id$ we have the commutative diagram with cartesian square
	\[
\begin{displaystyle}
   \xymatrix
   {	
	& \hat{X}_2\ar[d]\ar[ld] \ar[r]& X_2 \ar[d]\\
	\set{\Id} & Z_{\Bt, \set{\pm \Id}} \ar[l]^{s} \ar[r]_{\;\;\;\;r} & \Bt 
   }
\end{displaystyle}   
\]
where $\hat{X}_2 = r|_{Z_{\Bt, \set{\pm \Id}}}^*X_2 = \left\{(A, B, C, w, t) \in Z_{\Bt, \set{\pm \Id}} \times \CC^* \times \Bt\,|\, w^2 = t-2, t = \tr A \right\}$. Using that $Z_{\Bt, \set{\pm \Id}} \cong \Ss{\pm t_0} \times \SL{2}(\CC)$ we find that $\hat{X}_2 = \left\{w = \sqrt{\pm t_0 -2}\right\} \times \Ss{\pm t_0} \times \SL{2}(\CC)$ so $[\hat{X}_2] = 2(q^2+q)(q^3-q)$. This implies that $s_!r|_{Z_{\Bt, \set{\pm \Id}}}^*[X_2] = 2(q^2+q)(q^3-q)T_{\pm 2}$

	\item Over $[J_+]$, we a commutative diagram
	\[
\begin{displaystyle}
   \xymatrix
   {	
	& \hat{X}_2\ar[d]\ar[ld] \ar[r]& X_2 \ar[d]\\
	[J_+] & Z_{\Bt, [J_+]}^0 \ar[l]^{s} \ar[r]_{\;\;\;\;r} & \Bt 
   }
\end{displaystyle}   
\]
where $\hat{X}_2 = \left\{(A, w) \in Z_{\Bt, [J_+]}^0 \times \CC^* \,|\, w^2 = \tr A-2 \right\}$. Stratifying according to $z$ as in Lemma \ref{lem:fibr-w} we have the following strata.
\begin{itemize}
	\item For $z=0$, the equation for $W$ in Lemma \ref{lem:fibr-w} forces that $t = t_0$. In this way, we have that
\begin{align*}
	\hat{X}_2 \cap \set{z=0} &= \left\{(y, w) \in \left(Z_{\Bt, [J_+]}^0 \cap \set{z=0}\right) \times \CC^* \,\left|\, \begin{matrix}w^2 = t_0-2\\ y \neq 0 \end{matrix}\right.\right\} \\
	& = \CC^* \times \set{w = \pm \sqrt{t_0 -2}}.
\end{align*}
Therefore, $[\hat{X}_2 \cap \set{z=0}] = 2(q-1)$.
	\item For $z \neq 0$, we can solve for $y$ to obtain that
\begin{align*}
	\hat{X}_2 \cap \set{z \neq 0} &= \left\{(t, z, w) \in \left(\CC-\set{\pm 2}\right) \times \CC^* \times \CC^* \,\left|\, \begin{matrix}w^2 =  t-2\end{matrix}\right.\right\} \\
	&= X_2 \times \CC^*.
\end{align*}
The first factor is an affine parabola without $3$ points so we have that $[\hat{X}_2 \cap \set{z \neq 0}] = (q-1)(q-3)$.
\end{itemize}
Adding up all the contributions, we find that $[\hat{X}_2] = (q-1)^2$ and, thus, $s_!r|_{{Z_{\Bt, [J_+]}}}^*[X_2] = (q^2+q)(q^3-q)(q-1)^2T_+$.
	\item Over $[J_-]$, the computation is completely analogous and we have $s_!r|_{{Z_{\Bt, [J_-]}}}^*[X_2] = (q^2+q)(q^3-q)(q-1)^2T_-$.
	\item Over $\Bt$, the commutative diagram
	\[
\begin{displaystyle}
   \xymatrix
   {	
	& \hat{X}_2\ar[d]\ar[ld] \ar[r]& X_2 \ar[d]\\
	\Bt & Z_{\Bt, \Bt}^0 \ar[l]^{s} \ar[r]_{\;\;\;\;r} & \Bt 
   }
\end{displaystyle}   
\]
is completed with
$$
	\hat{X}_2 = \left\{(y, z, t,t', w) \in W \times \CC^*\,\left|\,\begin{matrix}w^2=t-2\\t,t' \neq \pm 2\end{matrix}\right.\right\}.
$$
The projection becomes $s(y,z,t,t', w)=t'$. As above, we stratify depending on $z$.

\begin{itemize}
	\item Where $z \neq 0$, the situation is simple since we can solve for $y$ and thus
	$$
	\hat{X}_2 \cap \set{z \neq 0} = \left\{\begin{matrix}w^2=t-2\\t,t' \neq \pm 2\end{matrix}\right\} \times \CC^* = \left\{\begin{matrix}w^2 = t-2\\ t \neq \pm 2\end{matrix}\right\} \times \CC^* \times \Bt.
$$
Moreover, under this isomorphism the map $s$ becomes the projection onto the last component. Therefore, since $\left\{w^2 = t-2, t \neq \pm 2\right\} \cong \CC - \set{0,\pm 2i}$, we obtain that $\left[\hat{X}_2 \cap \set{z \neq 0}\right] = (q-3)(q-1)T_{\Bt}$.
	\item Where $z = 0$, by Lemma \ref{lem:fibr-w} we get that this stratum contributes as $qT_{\Bt} + qS_2 + qS_{-2} + qS_{2}\times S_{-2} - 2qT_{t_0} - 2qT_{-t_0}$.
\end{itemize}

Adding up all the contributions and subtracting $\Zg{}(L_{\Ss{t_0}})(T_{\Bt})$, we finally find that
$$
	\Zg{}(L_{\Ss{t_0}})(S_2) = (q^2+q)(q^3-q)\left(T_2 + T_{-2} + (1-q)T_{\Bt} + qS_2 + qS_{-2} - qT_{t_0} - qT_{t_0}\right).
$$

\end{itemize}

\subsection{Image of $S_{-2}$} The calculations for this element are analogous for the ones of $S_2$, so we also obtain
$$
	\Zg{}(L_{\Ss{t_0}})(S_{-2}) = (q^2+q)(q^3-q)\left(T_2 + T_{-2} + (1-q)T_{\Bt} + qS_2 + qS_{-2} - qT_{t_0} - qT_{t_0}\right).
$$

\subsection{Image of $S_2 \times S_{-2}$} The computation for this element is very similar the one of $S_{\pm 2}$. As above, we will compute $s_!r^*[X_2 \times_{\Bt} X_{-2}]$ and, using that $[X_2 \times_{\Bt} X_{-2}] = T_{\Bt} + S_2 + S_{-2} + S_2 \times S_{-2}$, we get the image of $S_2 \times S_{-2}$ as $\Zg{}(L_{\Ss{t_0}})(S_2 \times S_{-2})=s_!r^*[X_2 \times_{\Bt} X_{-2}]-\Zg{}(L_{\Ss{t_0}})(T_{\Bt}) - \Zg{}(L_{\Ss{t_0}})(S_2) - \Zg{}(L_{\Ss{t_0}})(S_{-2})$.

As above, we decompose the calculation according to the stratification of $r^{-1}(\Bt)$.
\begin{itemize}
	\item Over $\pm \Id$ we have the commutative diagram with cartesian square
	\[
\begin{displaystyle}
   \xymatrix
   {	
	& \widehat{X_2 \times_{\Bt} X_{-2}}\ar[d]\ar[ld] \ar[r]& X_2 \times_{\Bt} X_{-2} \ar[d]\\
	\set{\Id} & Z_{\Bt, \set{\pm \Id}} \ar[l]^{s} \ar[r]_{\;\;\;\;r} & \Bt 
   }
\end{displaystyle}   
\]
where the pullback variety $\widehat{X_2 \times_{\Bt} X_{-2}} = X_2 \times_{\Bt} X_{-2} \times_{\Bt} Z_{\Bt, \set{\pm \Id}}$ is given by
\begin{align*}
	\widehat{X_2 \times_{\Bt} X_{-2}} &= \left\{(A, B, C, w_1, w_2, t) \in Z_{\Bt, \set{\pm \Id}} \times (\CC^*)^2 \times \Bt\,|\, w_1^2 = \tr A-2, w_2^2 = \tr A+2 \right\} \\
	& = \left\{w_1^2 = t_0-2, w_2^2 = t_0+2 \right\} \times \Ss{\pm t_0} \times \SL{2}(\CC).
\end{align*}
Therefore, we get $\left[\widehat{X_2 \times_{\Bt} X_{-2}}\right] = 4(q^2+q)(q^3-q)$.

\item Over $[J_+]$, there is a commutative diagram
	\[
\begin{displaystyle}
   \xymatrix
   {	
	& \widehat{X_2 \times_{\Bt} X_{-2}}\ar[d]\ar[ld] \ar[r]& X_2 \ar[d]\\
	[J_+] & Z_{\Bt, [J_+]}^0 \ar[l]^{s} \ar[r]_{\;\;\;\;r} & \Bt 
   }
\end{displaystyle}   
\]
where $\widehat{X_2 \times_{\Bt} X_{-2}} = \left\{(A, w_1, w_2) \in Z_{\Bt, [J_+]}^0 \times (\CC^*)^2 \,|\, w_1^2 = \tr A-2, w_2^2 = \tr A+2 \right\}$. Stratifying according to $z$, we get
\begin{itemize}
	\item For $z=0$, the equation for $H$ in Lemma \ref{lem:fibr-w} forces that $t = t_0$ and, thus
\begin{align*}
	\widehat{X_2 \times_{\Bt} X_{-2}} \cap \set{z=0} &= \left\{(y, w_1, w_2) \in \left(Z_{\Bt, [J_+]}^0 \cap \set{z=0}\right) \times (\CC^*)^2 \,\left|\, \begin{matrix}w_1^2 = t_0-2 \\ w_2^2 = t_0+2\\ y \neq 0 \end{matrix}\right.\right\} \\
	& = \CC^* \times \set{w_1 = \pm \sqrt{t_0 -2}, w_2 = \pm \sqrt{t_0 +2}}.
\end{align*}
Therefore, $\left[\widehat{X_2 \times_{\Bt} X_{-2}} \cap \set{z=0}\right] = 4(q-1)$.
	\item For $z \neq 0$, we solve for $y$ so that
\begin{align*}
	\widehat{X_2 \times_{\Bt} X_{-2}} \cap \set{z \neq 0} &= \left\{(t, z, w_1, w_2) \in \Bt \times (\CC^*)^3 \,\left|\, \begin{matrix}w_1^2 =  t-2\\ w_2^2 =  t+2 \end{matrix}\right.\right\} \\
	&= \left(X_2 \times_{\Bt} X_{-2} \right) \times \CC^*.
\end{align*}
The first factor is an affine hyperbola minus $4$ points, that contributes with $q-5$. Therefore, $\left[\widehat{X_2 \times_{\Bt} X_{-2}} \cap \set{z \neq 0}\right] = (q-1)(q-5)$.
\end{itemize}
Adding up all the contributions, we find that $\left[\widehat{X_2 \times_{\Bt} X_{-2}}\right] = (q-1)^2$ and, thus, $s_!r|_{{Z_{\Bt, [J_+]}}}^*\left[X_2 \times_{\Bt} X_{-2}\right] = (q^2+q)(q^3-q)(q-1)^2T_+$.
	\item Over $[J_-]$, the situation is completely analogous to the previous item so we have that $s_!r|_{{Z_{\Bt, [J_-]}}}^*\left[X_2 \times_{\Bt} X_{-2}\right] = (q^2+q)(q^3-q)(q-1)^2T_+$.
	\item Over $\Bt$, we have the commutative diagram
	\[
\begin{displaystyle}
   \xymatrix
   {	
	& \widehat{X_2 \times_{\Bt} X_{-2}}\ar[d]\ar[ld] \ar[r]& X_{2}\times X_{-2} \ar[d]\\
	\Bt & Z_{\Bt, \Bt}^0 \ar[l]^{s} \ar[r]_{\;\;\;\;r} & \Bt 
   }
\end{displaystyle}   
\]
that is completed with $
	\widehat{X_2 \times_{\Bt} X_{-2}} = \left\{w_1^2=t-2, w_2^2=t+2\right\}$. As always, we stratify depending on $z$.

\begin{itemize}
	\item Where $z \neq 0$, we get a trivial fibration as for $S_2$ with fiber $X_2 \times_{\Bt} X_{-2}$ so we obtain that $\left[\widehat{X_2 \times_{\Bt} X_{-2}} \cap \set{z \neq 0}\right] = (q-5)(q-1)T_{\Bt}$.
	\item Where $z = 0$, by Lemma \ref{lem:fibr-w} we get that this stratum contributes as $2qT_{\Bt} + 2qS_{2} + 2qS_{-2} + 2qS_{2} \times S_{-2} - 4qT_{t_0} - 4qT_{-t_0}$.
\end{itemize}
\end{itemize}

Therefore, putting together these results we find that
$$
	\Zg{}(L_{\Ss{t_0}})(S_2 \times S_{-2}) = (q^2+q)(q^3-q)\left(T_2 + T_{-2} + T_{\Bt} + qS_2 \times S_{-2} - qT_{t_0} - qT_{t_0}\right).
$$

\subsection{Image of $T_{t_0}$}\label{sec:img-t0}

\begin{itemize}
	\item Over $\Id$, we have that $Z_{\Ss{t_0}, \set{\Id}} = \Ss{t_0} \times \SL{2}(\CC)$, so we have that $\Zg{}(L_{\Ss{t_0}})(T_{t_0})|_{\Id} = (q^2+q)(q^3-q) T_2$.
	\item Over $-\Id$, we have that $Z_{\Ss{t_0}, \set{-\Id}} = \emptyset$, except in the case that $t_0 = 0$ for which $Z_{\Ss{t_0}, \set{-\Id}} = \Ss{t_0} \times \SL{2}(\CC)$ (c.f.\ Remark \ref{rmk:trace-equal-0}). Hence, $\Zg{}(L_{\Ss{t_0}})(T_{t_0})|_{-\Id} = 0$ if $t_0 \neq 0$ and $\Zg{}(L_{\Ss{0}})(T_{0})|_{-\Id} = (q^2+q)(q^3-q) T_{-2}$.
	\item Over $[J_+]$, the calculation is a simplified version of the one in Section \ref{sec:image-T-Bt}. On $W$, having $t = t_0$ and $t' = 2$ implies that $yz=0$. Hence, $Z_{\Ss{t_0}, [J_+]}^0 = \set{yz=0} - \set{(0,0)} = \CC^* \sqcup \CC^*$ so $\Zg{}(L_{\Ss{t_0}})(T_{t_0})|_{[J_+]} = (q^2+q)(q^3-q)(2q-2) T_+$.
	\item Over $[J_-]$, now we have that, in $W$, to have $t = t_0$ and $t' = -2$ implies that $yz = \alpha_0$ for some $\alpha_0 \in \CC$ fixed. If $t_0 \neq 0$, then $\alpha_0 \neq 0$. Thus $Z_{\Ss{t_0}, [J_-]}^0 = \set{yz=\alpha_0} = \CC^*$ so $\Zg{}(L_{\Ss{t_0}})(T_{t_0})|_{[J_-]} = (q^2+q)(q^3-q)(q-1) T_-$. However, if $t_0=0$ then $\alpha_0 = 0$ and we have the same situation as for $[J_+]$, so $\Zg{}(L_{\Ss{0}})(T_{0})|_{[J_-]} = (q^2+q)(q^3-q)(2q-2) T_-$.
	 	
	\item Over $\Bt$, we also have a simplified version of Section \ref{sec:image-T-Bt}. In this case, $Z_{\Ss{t_0}, \Bt}^0$ is given by the tuples $(t', y, z) \in \Bt \times \CC^2$ such that
$$
\frac{(t'-t_0)^2}{t_0^2-4} - \frac{t't_0}{t_0+2} +yz = -1.
$$
\begin{itemize}
	\item For $z=0$, this equation forces $t' = t_0^2-2$. Hence, this stratum contributes as $qT_{t_0^2-2}$ if $t_0 \neq 0$, and $0$ if $t_0=0$.
	\item For $z \neq 0$, we can solve the equation for $y$ and we get $Z_{\Ss{t_0}, \Bt}^0 \cap \set{z \neq 0} = \Bt \times \CC^*$, being the map $s$ the projection onto the first component. Thus, this stratum contributes with $(q-1) T_{\Bt}$. 
\end{itemize}
\end{itemize}

Adding up all the contributions, we get that, for $t_0 \neq 0$
$$
	\Zg{}(L_{\Ss{t_0}})(T_{t_0}) = (q^2+q)(q^3-q)\left(T_2 + (2q-2) T_+ + (q-1)T_- + (q-1)T_{\Bt} + qT_{t_0^2-2}\right).
$$

On the other hand, in the case $t_0 = 0$ we obtain
$$
	\Zg{}(L_{\Ss{0}})(T_{0}) = (q^2+q)(q^3-q)\left(T_2 + T_{-2} + (2q-2) T_+ + (2q-2)T_- + (q-1)T_{\Bt}\right).
$$

\subsection{Image of $T_{-t_0}$}\label{sec:img-mt0} The computation for this generator is analogous to the one of $T_{t_0}$. The differences are that, now, the stratum for $\Id$ is empty, except in the case $t_0 = 0$, and the fibration over $\Bt$ in the case $z=0$ lands in $2-t_0^2$. Taking into account these changes, we obtain that, for $t_0 \neq 0$
$$
	\Zg{}(L_{\Ss{t_0}})(T_{-t_0}) = (q^2+q)(q^3-q)\left(T_{-2} + (q-1) T_+ + (2q-2)T_- + (q-1)T_{\Bt} + qT_{2-t_0^2}\right).
$$

\begin{rmk}
In the case $t_0=0$, the generators $T_{t_0}$ and $T_{-t_0}$ are indistinguishable, so both calculations of Sections \ref{sec:img-t0} and \ref{sec:img-mt0} must apply. For this reason, $\Zg{}(L_{\Ss{0}})(T_{0})$ has to be symmetric in $T_2, T_{-2}$ and $T_+, T_-$. This explains the balanced shape of $\Zg{}(L_{\Ss{0}})(T_{0})$.
\end{rmk}

\subsection{Image of $T_{t}$ with $t \neq \pm t_0$}
\begin{itemize}
	\item Over $\pm \Id$, we have that $Z_{\Ss{t}, \set{\pm \Id}} = \emptyset$, so these strata add no contribution.
	\item Over $J_\pm$, we have that $Z_{\Ss{t}, [J_\pm]}^0 = W \cap \set{t'=2, t\textrm{ fixed}}$. In particular, since $t \neq \pm t_0$, the product $yz \neq 0$ so we can solve for $y$. Hence, we obtain that $Z_{\Ss{t_0}, [J_\pm]}^0 = \CC^*$ so these strata contributes with $(q-1)(q^2+q)(q^3-q)T_{\pm}$.
	\item Over $\Bt$, we have that $Z_{\Ss{t}, \Bt}^0 = W \cap \set{t\textrm{ fixed}}$.
\begin{itemize}
	\item For $z=0$, $t'$ is forced to be one of the two (different) roots of the equation for $W$ in Lemma \ref{lem:fibr-w}, namely $\tau_+(t_0,t), \tau_-(t_0,t) \in \CC$. Hence, this stratum contributes as $qT_{\tau_+(t_0, t)} + qT_{\tau_-(t_0, t)}$.
	\item For $z \neq 0$, we can solve this equation for $y$ and we get $Z_{\Ss{t}, \Bt}^0 \cap \set{z \neq 0} = \Bt \times \CC^*$, being the map $s$ the projection onto the first component. Thus, this stratum contributes with $(q-1) T_{\Bt}$. 
\end{itemize}

\end{itemize}

Therefore, taking into account all the contributions, we finally get that
\begin{align*}
	\Zg{}(L_{\Ss{t_0}})(T_{t}) =& (q^2+q)(q^3-q)\left((q-1) T_+ + (q-1)T_- + (q-1)T_{\Bt}+ qT_{\tau_+(t_0,t)} + qT_{\tau_-(t_0,t)} \right).
\end{align*}

\begin{rmk}
The new traces $\tau_{\pm}(t_0, t)$ can be better understood as follows. Let us consider the eigenvalues of $t$ and write it as $t=\lambda + \lambda^{-1}$, for $\lambda \in \CC^*-\set{\pm 1}$. Then, the two new traces $\tau_{\pm}(t_0, t)$ are the ones associated to the eigenvalues $\lambda_0\lambda$ and $\lambda_0\lambda^{-1}$, that is, $\tau_+(t_0, t)=\lambda_0\lambda+\lambda_0^{-1}\lambda^{-1}$ and $\tau_-(t_0, t)=\lambda_0\lambda^{-1}+\lambda_0^{-1}\lambda$ (or vice-versa, the names $\tau_{\pm}$ are arbitrary).
In particular, for $t = \pm t_0$ we have that $\tau_+(t_0, \pm t_0) = \pm (t_0^2-2)$ and $\tau_-(t_0,\pm t_0) = \pm 2$. Since we have to dismiss the root $\pm 2$, we only have one contribution in Sections \ref{sec:img-t0} and \ref{sec:img-mt0}.
\end{rmk}

\subsection{The reduced TQFT} From the computations of the previous sections, we have obtained an explicit expression of the morphism $\Zg{}(L_{\Ss{t_0}})$ on $\overline{\cW}$. With this information at hand, passing to the localization with respect to $q, q+1$ and $q-1$ as mentioned in Remark \ref{rmk:loc-kvar}, we can finally compute the reduced TQFT, $\cZg{}(L_{\Ss{t_0}})$.

\begin{thm}\label{thm:image-Lt}
Fix $t_0 \in \Bt=\CC-\set{\pm 2}$. Under the module homomorphism $\cZg{}(L_{\Ss{t_0}}): \KoVarrel{[\SL{2}(\CC)/\SL{2}(\CC)]} \to \KoVarrel{[\SL{2}(\CC)/\SL{2}(\CC)]}$, we have that $\cZg{}(L_{\Ss{t_0}})\left(\overline{\cW}\right) \subseteq \overline{\cW}$. Explicitly, in the set of generators, the matrix is
\tiny
$$
(q^3-q)\left(\begin{array}{c|ccccccccccc}
&T_2 & T_{-2} & T_+ & T_- & T_{\Bt} & S_2 & S_{-2} & S_2 \times S_{-2} & T_{t_0} &
T_{-t_0} & T_{t} \\\hline
T_2 & 0 & 0 & 0 & 0 & 1 & 1 & 1 & 1 & 1 &
0 & 0  \\
T_{-2} &0 & 0 & 0 & 0 & 1 & 1 & 1 & 1 & 0 &
1 & 0  \\
T_+ & 0 & 0 & q & q & q^2-q & 0
& 0 & 1-q & 2 \, q - 2 &  q - 1  & q - 1 \\
T_- & 0 & 0 & q & q & q^2-q & 0
& 0 & 1-q & q - 1 & 2 \, q - 2  & q - 1 \\
T_{\Bt} & 0 & 0 & q & q & q^{2} - q + 1 & 1-q &
1-q & 2-q & q - 1 & q - 1 & q - 1 \\
S_2 & 0 & 0 & 0 & 0 & 0 & q & q & 0 & 0 &
0 & 0   \\
S_{-2} &0 & 0 & 0 & 0 & 0 & q & q & 0 & 0 &
0 & 0   \\
S_2 \times S_{-2} &0 & 0 & 0 & 0 & q & 0 & 0 & q & 0 &
0 & 0   \\
T_{t_0} & q^2+q & 0 & q & 0 & -q & -q & -q & -q &
0 & 0 & 0   \\
T_{-t_0} &0 & q^2+q & 0 & q & -q & -q & -q & -q &
0 & 0 & 0 \\
T_{t_0^2-2} & 0 & 0 & 0 & 0 & 0 & 0 & 0 & 0 & q &
0  & 0  \\
T_{2-t_0^2} &0 & 0 & 0 & 0 & 0 & 0 & 0 & 0 & 0 &
q  & 0 \\
T_{\tau_{+}(t_0,t)} & 0 & 0 & 0 & 0 & 0 & 0 & 0 & 0 & 0 &
0  & q \\
T_{\tau_{-}(t_0,t)} & 0 & 0 & 0 & 0 & 0 & 0 & 0 & 0 & 0 &
0  & q 
\end{array}\right)
$$
\normalsize
In the case $t_0 = 0$, the images of $T_{\pm 2}, T_\pm, S_{\pm 2}, S_2 \times S_{-2}$ or $T_t$ with $t \neq 0$ are the same as above, but the image of $T_0$ is given by
$$
	\cZg{}(L_{\Ss{t_0}})(T_{0}) = (q^3-q)\left(T_2 + T_{-2} + (2q-2) T_+ + (2q-2)T_- + (q-1)T_{\Bt}\right)
$$

\begin{proof}
From the computations above, we know the explicit expression of $\Zg{}(L_{\Ss{t_0}}): \overline{\cW} \to \overline{\cW}$. The reduced TQFT is $\cZg{}(L_{\Ss{t_0}}) = \Zg{}(L_{\Ss{t_0}}) \circ \eta^{-1}$. Thus, as $\eta(\overline{\cW}) \subseteq \overline{\cW}$, the same holds for $\cZg{}(L_{\Ss{t_0}})$. Moreover, from the expression of $\eta$ given by \cite[Proposition 6.3]{GP:2018a} and Lemma \ref{lem:eta-morph}, the matrix of $\cZg{}(L_{\Ss{t_0}})$ can be explicitly computed using a computer algebra system.
\end{proof}
\end{thm}

\subsection{Image of $T_t$ under other tubes}
\label{sec:image-Tt-other-tubes}

Once completed the computation of $\Zg{}(L_{\Ss{t_0}})$, we need to extend the calculations in Sections 6.3 and 6.4 of \cite{GP:2018a} for $\Zg{}(L)$ and $\Zg{}(L_{[J_\pm]})$ to the whole $\overline{\cW}$. Since the image of the generators of $\cW \subseteq \overline{\cW}$ is known \cite[Theorem 6.5]{GP:2018a}, it is enough to compute the elements $\Zg{}(L)(T_{{t}})$, $\Zg{}(L_{[J_+]})(T_{t})$ and $\Zg{}(L_{[J_-]})(T_{{t}})$.

This task can be performed following the strategy of the previous sections, but there is a shortest path that we shall explain here. By Theorem \ref{thm:image-Lt}, we have $\cZg{}(L_{\Ss{t}})(T_2)=(q^2+q)(q^3-q)T_t$.
Now, observe that the bordisms $L_{[J_\pm]}: (S^1, \star) \to (S^1, \star)$ and $L_{\Ss{t}}: (S^1, \star) \to (S^1, \star)$ commute as morphisms of $\Bordppar{2}{\Lambda}$ so $\cZg{}(L_{\Ss{t}}) \circ \cZg{}(L_{[J_\pm]}) = \cZg{}(L_{[J_\pm]}) \circ \cZg{}(L_{\Ss{t}})$. Therefore, since $\eta(T_t)=(q^2+q)T_t$ and $\cZg{}(L_{\Ss{t}})(T_2)=T_t$, we get that
\begin{align*}
	(q^3-q)\Zg{}(L_{[J_\pm]})(T_t) &= (q^3-q)\cZg{}(L_{[J_\pm]})\left(\eta(T_t)\right) = (q^3-q)(q^2+q)\cZg{}(L_{[J_\pm]})\left(T_t\right) \\
	& = \cZg{}(L_{[J_\pm]}) \circ \cZg{}(L_{\Ss{t}})(T_2) =  \cZg{}(L_{\Ss{t}}) \circ \cZg{}(L_{[J_\pm]})(T_2).
\end{align*}
These last term can be computed from the expression of $\cZg{}(L_{[J_\pm]})$ in \cite[Theorem 6.5]{GP:2018a} and Theorem \ref{thm:image-Lt} for $\cZg{}(L_{\Ss{t}})$. Moreover, analogous arguments hold for the bordism $L: (S^1, \star) \to (S^1, \star)$. Putting together this information we obtain the following result.

\begin{prop} The image of $T_t \in \KVarrel{[\SL{2}(\CC)/\SL{2}(\CC)]}$ under the endomorphisms $\Zg{}(L_{[J_\pm]})$ and $\Zg{}(L)$ is
\begin{align*}
	\Zg{}(L_{[J_\pm]})(T_t) =&\, (q^3-q)^2 \left(T_+ + T_- + T_{\Bt} + T_{\pm t}\right), \\
	\Zg{}(L)(T_t)  = &\, (q^3-q)^2\left((q^2+4q+1) (T_2 + T_{-2})+ (q^2+2q+3)(q-1)q (T_+ + T_-) \right.\\
	&+ (q^4+q^3-q^2+q+1)T_{\Bt} +3q^2 (S_2 + S_{-2}) + (q^3+q^2+q) S_2 \times S_{-2} \\
	& + \left.(q^3-q^2)T_{t}\right). 
\end{align*}
\end{prop}

\begin{rmk}
As in the case of Theorem \ref{thm:image-Lt}, from this information and the explicit expression of the endomorphism $\eta$, we can easily compute the reduced versions $\cZg{}(L)$ and $\cZg{}(L_{[J_\pm]})$.
\end{rmk}

\section{The interaction phenomenon}
\label{sec:interation-phenomenon}

The previous computation shows that there exists a dichotomy in the behavior of the TQFT on the skyscraper generators over $\Bt$. Fix $t_0 \in \Bt$ and consider the tube $L_{\Ss{t_0}}: (S^1, \star) \to (S^1, \star)$. As shown in Theorem \ref{thm:image-Lt}, for $t \in \Bt - \set{ \pm t_0}$, the images of the generators $T_t$ under $L_{\Ss{t_0}}$ are
\begin{align*}
	\cZg{}(L_{\Ss{t_0}})(T_{t}) =\, &(q^2+q)(q^3-q)\big((q-1)T_+ + (q-1)T_- + (q-1)T_{\Bt} 
	\\& + qT_{\tau_+(t_0, t)}+ qT_{\tau_-(t_0, t)}\big).
\end{align*}
Therefore, these images are essentially equal with varying $t \in \Bt - \set{\pm t_0}$. However, for $t = t_0$ the image has a different shape 
\begin{align*}
	\cZg{}(L_{\Ss{t_0}})(T_{t_0}) =\, & (q^2+q)(q^3-q)\big( (q-1)T_{+} + (q-1)T_{-} + (q-1)T_{\Bt} \\
	& + qT_{t_0^2-2}+T_{2} + (q-1)T_{+}\big).
\end{align*}
We can interpret this last element as follows. In this case, we have that $\tau_+(t_0, t_0) = t_0^2-2$ but $\tau_-(t_0, t_0) = 2$. Hence, one of the two new skyscraper generators lies over the matrices of trace $2$ i.e.\ $\set{\Id} \cup [J_+]$. In this way, the contribution $qT_{\tau_-(t_0, t_0)}$ that appears in the generic case becomes $T_{2} + (q-1)T_{+}$ in the case $t = t_0$. We may understand this phenomenon as that the trace $t_0$ interacts with the generator $T_{t_0}$ destroying the element $qT_{\tau_-(t_0, t_0)}$ and creating two new elements $T_{2}$ and $(q-1)T_{+}$. Analogous considerations can be done for $t = -t_0$, where $qT_{\tau_-(t_0, -t_0)}$ is destroyed and a new contribution $T_{-2} + (q-1)T_{-}$ appears.

This `interaction phenomenon' makes the computation of the iterated element $\cZg{}(L_{\Ss{t_s}}) \circ \ldots \circ \cZg{}(L_{\Ss{t_1}})(T_2) \in \overline{\cW}$ much more involved that the non-semisimple case of \cite{GP:2018a}. In this section, we shall study thoroughly this interaction phenomenon and we shall provide an explicit formula for the iterated element.

\subsection{The generic part}

Let us consider in the first place the generic case. The importance of this situation is that, for generic traces $t_1, \ldots, t_s$, we can obtain a closed formula for the image of the bordism $L_{\Ss{t_s}} \circ \ldots \circ L_{\Ss{t_1}} \circ D: \emptyset \to (S^1, \star)$.

\begin{defn}\label{defn:generic-traces}
Let $t_1, \ldots, t_s \in \Bt$ be traces, and let us write them as $t_i = \lambda_i + \lambda_i^{-1}$ for some $\lambda_i \in \CC^*-\set{\pm 1}$. We will say that they are \emph{generic} if all the products $\lambda_1^{\pm 1}\cdots\lambda_\ell^{\pm 1} \neq \pm 1$ for all $1 \leq \ell \leq s$.
\end{defn}

\begin{rmk}
This definition of genericity agrees with the Definition 4.6.1 of \cite{Mellit:2019}. However, as we will see, we came up to this definition guided through a different path.
\end{rmk}

With a view towards the following result, for simplicity let us denote by $\tau: \CC^*-\set{\pm 1} \to \Bt=\CC-\set{\pm 2}$ the map $\tau(\lambda)=\lambda + \lambda^{-1}$. This is the quotient map for the action of $\ZZ_2$ on $\CC^*-\set{\pm 1}$ by $\lambda \mapsto \lambda^{-1}$, i.e.\ permutation of the eigenvalues.

\begin{prop}\label{prop:generic-case}
Let $t_1, \ldots, t_s \in \Bt$ be a generic set of traces, and let us write them as $t_i = \lambda_i + \lambda_i^{-1}$ for some $\lambda_i \in \CC^*-\set{\pm 1}$. Then we have
\begin{align*}
	\cZg{}(L_{\Ss{t_s}} \circ \ldots \circ L_{\Ss{t_1}})(T_{2}) =\, & (q^3-q)^{s-1}\left(a_s (T_2 + T_{-2}) + b_s (T_+ + T_-) + c_s T_{\Bt} + d_s S_2 \times S_{-2}\right)\\
	&+ \frac{(q^3-q)^s(q+1)}{2} \sum_{(\epsilon_1, \ldots, \epsilon_s)  \in \set{\pm 1}^s} q^sT_{\tau(\lambda_1^{\epsilon_1}\cdots \lambda_s^{\epsilon_s})}.
\end{align*}
The coefficients are given by
\begin{align*}
	a_s & = q^s(q + 1)^s + q^2(q + 1)^s - 2^{s - 1}q^s(q + 1)^2,\\
	b_s & = q^s(q + 1)^{s+1}(q-1) - 2^{s-1}q^s(q + 1)^2(q-1),\\
	c_s & = q^{s+2}(q + 1)^s + q^2(q + 1)^s-2^{s-1}q^{s+1}(q + 1)^2,\\
	d_s & = q^{s+1}(q + 1)^s + q^3(q + 1)^s - 2^{s-1}q^{s+1}(q+ 1)^2.
\end{align*}
\begin{proof}
The formulae for the coefficients $a_s, b_s, c_s$ and $d_s$ follow easily by an induction argument using the matrix of Theorem \ref{thm:image-Lt}. For the coefficients of the skyscraper generators $T_t$ for $t \in \Bt$, the key point is to observe that the linear relation $(q^2+q)a_s + qb_s - qc_s -q d_s=0$ holds for any $s > 0$. Hence, working by induction we have that
\begin{align*}
\cZg{}(L_{\Ss{t_{s+1}}} &\circ \ldots \circ L_{\Ss{t_1}})(T_{2}) = (q^3-q)^{s}\left(a_{s+1} (T_2 + T_{-2})  + b_{s+1} (T_+ + T_-) \right.\\& \left. + c_{s+1} T_{\Bt} + d_{s+1} S_2 \times S_{-2} + \left((q^2+q)a_s + qb_s - qc_s -q d_s\right) T_{t_{s+1}}\right)
	\\ &+ \frac{(q^3-q)^{s+1}(q+1)}{2} \sum_{\epsilon_1, \ldots, \epsilon_s} q^s\left.\cZg{}(L_{\Ss{t_{s+1}}})\left(T_{\tau(\lambda_1^{\epsilon_1}\cdots \lambda_s^{\epsilon_s})}\right)\right|_{\Bt}.
\end{align*}
The formula for the skyscraper generators follows from the observation that, for any $\lambda \in \CC^* - \set{\pm 1}$, we have $\left.\cZg{}(L_{\Ss{t_{s+1}}})\left(T_{\tau(\lambda)}\right)\right|_{\Bt} = q T_{\tau(\lambda\lambda_{s+1})} + qT_{\tau(\lambda\lambda_{s+1}^{-1})}$.
\end{proof}
\end{prop}

\subsection{The interaction term}
\label{sec:interaction-term}

Despite our success in the computation of the generic situation, we still need to correct this result to take into account the contribution in the non-generic case. For this purpose, suppose that we have traces $t_1, \ldots, t_s \in \Bt$ and let us write them as $t_i = \lambda_i + \lambda_i^{-1}$ for $\lambda_i \in \CC^*-\set{\pm 1}$. We shall denote
\begin{align*}
	\alpha_+ &= \frac{1}{2} \left|\left\{(\epsilon_1, \ldots, \epsilon_s) \in \set{\pm 1}^s \,\left|\, \begin{matrix}\lambda_1^{\epsilon_1}\cdots \lambda_s^{\epsilon_s} = 1 \end{matrix} \right. \right\}\right|,\\ 
\alpha_- &= \frac{1}{2} \left|\left\{(\epsilon_1, \ldots, \epsilon_s) \in \set{\pm 1}^s \,\left|\, \begin{matrix}\lambda_1^{\epsilon_1}\cdots \lambda_s^{\epsilon_s} = -1 \end{matrix} \right. \right\}\right|.
\end{align*}
If we need to make explicit the set of traces over which $\alpha_\pm$ is computed, we will write $\alpha_\pm(t_1, \ldots, t_s)$.

In this situation, we define the \emph{interaction term} of the traces $t_1, \ldots, t_s$ as
$$
	\cI(t_1, \ldots, t_s) = q^{s-1}(q^3-q)^s(q+1)\left(\alpha_+(T_2 + (q-1)T_+ )+ \alpha_- (T_{-2} + (q-1) T_-)\right).
$$
Observe that $\cI$ does not depend on the chosen ordering of the traces $t_i$ nor the choice the eigenvalue $\lambda_i$ or $\lambda_i^{-1}$ for $t_i$. Hence, indeed it only depends on the set of traces.

The importance of this term comes from the following result.
\begin{prop}\label{prop:non-generic}
Let $t_1, \ldots, t_s \in \Bt$ be any traces. Then we have
\begin{align*}
\cZg{}(L_{\Ss{t_s}} \circ \ldots \circ L_{\Ss{t_1}})(T_{2}) =\, &(q^3-q)^{s-1}\left(a_s (T_2 + T_{-2}) + b_s (T_+ + T_-) + c_s T_{\Bt} + d_s S_2 \times S_{-2}\right) \\&
	+ \frac{(q^3-q)^s(q+1)}{2} \sum_{\epsilon_1, \ldots, \epsilon_s} q^sT_{\tau(\lambda_1^{\epsilon_1}\cdots \lambda_s^{\epsilon_s})}+ \cI(t_1, \ldots, t_s),
\end{align*}
where the sum on the skyscraper generators runs over the set of tuples $(\epsilon_1, \ldots, \epsilon_s) \in \set{\pm 1}^s$ such that $\lambda_1^{\epsilon_1} \cdots \lambda_s^{\epsilon_s} \neq \pm 1$.
\begin{proof}
Firstly, observe that, if $t_1, \ldots, t_s$ are generic, then $\cI(t_1, \ldots, t_s)=0$ and the result follows from Proposition \ref{prop:generic-case}. We will prove the result by induction on $s$. For the base case $s=1$, the formula holds since this case is always generic.

In the general case, we may suppose that $t_{s+1}$ is not generic with $t_{1}, \ldots, t_{s}$. For simplicity, let us denote $A_s = (q^3-q)^{s-1}(a_s (T_2 + T_{-2}) + b_s (T_+ + T_-) + c_s T_{\Bt} + d_s S_2 \times S_{-2})$ and $S_s = \frac{1}{2}\sum_{} q^sT_{\tau(\lambda_1^{\epsilon_1}\cdots \lambda_s^{\epsilon_s})}$. Observe that, if we shorten $\alpha_\pm^s = \alpha_\pm(t_1, \ldots, t_s)$ and $\alpha_\pm^{s+1} = \alpha_\pm(t_1, \ldots, t_{s+1})$, then the number of times that $T_{\pm t_0}$ appears in $S_s$ is precisely $\alpha_\pm^{s+1}$. Now we separate the sum as
$$
	S_s =  \alpha_+^{s+1} q^s T_{t_{s+1}} + \alpha_-^{s+1} q^s T_{-t_{s+1}} + \frac{1}{2}\sum_{\lambda_1^{\pm 1}\cdots \lambda_s^{\pm 1} \neq \lambda_{s+1}^{\pm 1}} q^sT_{\tau(\lambda_1^{\pm 1}\cdots \lambda_s^{\pm 1})}.
$$
When we compare $\cZg{}(L_{\Ss{t_s}})(S_s)$ with respect to the image of the same sum in the generic case, we see that each generator $T_{\pm t_{s+1}}$ gives an extra contribution to the submodule $\cW \subseteq \overline{\cW}$ of $(q^3-q)(T_{\pm 2} + (q-1)T_{\pm})$. This implies a total extra contribution of $q^{s}(q^3-q)^{s+1}(q+1)(\alpha_+^{s+1} (T_{2} + (q-1)T_{+}) + \alpha_-^{s+1} (T_{-2} + (q-1)T_{-})) = \cI(t_1, \ldots, t_{s+1})$ with respect to the generic case. This justifies the appearance of the new interaction term.

On the other hand, $S_s$ has $2^s-2\alpha_+^s -2\alpha_-^s$ summands, which is $2(\alpha_+^s + \alpha_-^s)$ fewer terms than the generic case. This implies that $\cZg{}(L_{\Ss{t_s}})(S_s)$ has a missing contribution of $(\alpha_+^s + \alpha_-^s)q^s(q^3-q)(q-1)(T_+ + T_- + T_{\Bt})$ to the submodule $\cW \subseteq \overline{\cW}$ with respect to the generic case. Moreover, it also creates a deficit of $2q^{s+1}(q^3-q)^{s+1}(q+1)(\alpha_+^s T_{t_{s+1}} + \alpha_-^sT_{t_{s+1}})$ to the submodule generated by the skyscraper generators, coming from the $2(\alpha_+^s + \alpha_-^s)$ missing terms in $S_s$ satisfying $\lambda_1^{\epsilon_1} \cdots \lambda_s^{\epsilon_s} = \pm 1$.
This missing contribution is offset with the image of the old interaction term, as
\begin{align*}
	\cZg{}(L_{\Ss{t_{s+1}}})\left(\cI(t_1, \ldots, t_s)\right) =\, &(\alpha_+^s + \alpha_-^s)q^{s}(q^3-q)^{s+1}(q^2-1)(T_+ + T_- + T_{\Bt}) \\
	& + 2q^{s+1}(q^3-q)^{s+1}(q+1)(\alpha_+^s T_{t_{s+1}} + \alpha_-^sT_{t_{s+1}}).
\end{align*}
Putting together these observations, we finally get that
\begin{align*}
	\cZg{}(L_{\Ss{t_{s+1}}})&\left(A_s + (q^3-q)^s(q+1)S_s + \cI(t_1, \ldots, t_s)\right) = A_{s+1} + (q^3-q)^{s+1}(q+1)S_{s+1}
	\\&+ \cI(t_1, \ldots, t_{s+1}) -\cZg{}(L_{\Ss{t_{s+1}}})(\cI(t_1, \ldots, t_s)) + \cZg{}(L_{\Ss{t_{s+1}}})(\cI(t_1, \ldots, t_s)) \\
	 & =A_{s+1}+ (q^3-q)^{s+1}(q+1)S_{s+1}+ \cI(t_1, \ldots, t_{s+1}),
\end{align*}
as we wanted to prove.
\end{proof}
\end{prop}

\begin{rmk}
In Definition \ref{defn:generic-traces}, we defined a set of traces $t_1, \ldots, t_s$ to be generic as the property of having no interaction within any subset $t_1, \ldots, t_\ell$ for $1 \leq \ell \leq s$. However, Proposition \ref{prop:non-generic} proves that, indeed, the final result does not depend on the potential interactions on proper subsets but only on the whole set $t_1, \ldots, t_s$. This is a quite surprising memoryless property that cannot be expected from the only shape of the TQFT.
\end{rmk}

Once we know how to control the interaction phenomenon, we can finally give a closed formula for the virtual class of any parabolic representation variety.

\begin{thm}\label{thm:result-complete} Let $\Sigma_g$ be the closed orientable genus $g$ surface and fix traces $t_1, \ldots, t_s$, maybe non-generic. Let $Q$ be a parabolic structure with $r$ punctures with holonomy $[J_+]$ and $s>0$ punctures with holonomies $\Ss{t_1}, \ldots, \Ss{t_s}$. The virtual class of $\Rep{\SL{2}(\CC)}(\Sigma_g, Q)$ in $\Ko{\Var{\CC}}$ is
\begin{itemize}
	\item If $r>0$, then
\begin{align*}
	\left[\Rep{\SL{2}(\CC)}(\Sigma_g, Q)\right] =\, & q^{2g + s-1}(q - 1)^{2g + r-1}(q+1)\left(2^{2g+s-1} - 2^s + (q + 1)^{2g + r+s-2} \right)\\
	&+ \overline{\cI}_r(t_1, \ldots, t_s).
\end{align*}
The interaction term is given by
\begin{align*}
\overline{\cI}_r(t_1, \ldots, t_s) =\, & q^{2g + s-1}(q - 1)^{2g + r-1}(\alpha_+ + \alpha_-)\bigg( 2^{2g} + 2^{2g}q  - 2q -2 \\
&  + (q + 1)^{2g + r} + (q+1)\left(1 - 2^{2g-1} - \frac{1}{2}(q + 1)^{2g + r-1}\right)\bigg) \\
& + q^{2g + s-1}(q - 1)^{2g + r}(q+1)\alpha_+.
\end{align*}
	\item If $r=0$, then
\begin{align*}
	\left[\Rep{\SL{2}(\CC)}(\Sigma_g, Q)\right] =\, & q^{2g+s-1}(q - 1)^{2g-1}(q + 1)(2^{2g+ s - 1} - 2^s + (q + 1)^{2g+s-2} \\
	&+ q^{2-2g-s}(q + 1)^{2g+s-2}) + \overline{\cI}_0(t_1, \ldots, t_s).
\end{align*}
The interaction term is given by
\begin{align*}
\overline{\cI}_0(t_1, \ldots, t_s) =\, & q^{s-1}(q - 1)^{2g-1}(q+1)(\alpha_+ + \alpha_-)\left(q(q + 1)^{2g-1} + q^{2g}(q + 1)^{2g-1}\right.
\\ 
&\left.- q^{2g}(q + 1)^{2g-1} - q(q + 1)^{2g-1}\right) +q^{2g+s-1}(q - 1)^{2g}(q + 1)\alpha_+.
\end{align*}
\end{itemize}
\begin{proof}
The proof is a computation using the results of \cite[Theorem 6.5 and Corollary 6.9]{GP:2018a} and Proposition \ref{prop:non-generic}. Let us focus on the case $r>0$, being the case $r=0$ analogous. In that case, if $A \subseteq \Sigma_g$ is a set of $g+r+s+1$ base points, we can decompose $(\Sigma_g, A, Q) = D^\dag \circ L^g \circ L_{[J_+]}^r \circ L_{\Ss{t_s}} \circ \ldots \circ L_{\Ss{t_1}} \circ D$. Therefore, since $\cZg{}(D)(1)=T_2$, via the reduced TQFT we obtain that
$$
	\left[\Rep{}(\Sigma_g, Q)\right] = \frac{1}{(q^3-q)^{g+r+s}} \cZg{}(D^\dag) \circ \cZg{}(L)^g \circ \cZg{}(L_{[J_+]})^r \circ \cZg{}(L_{\Ss{t_{s}}} \circ \ldots \circ L_{\Ss{t_{1}}})(T_2).
$$ 

The computation of $\cZg{}(L_{\Ss{t_{s}}} \circ \ldots \circ L_{\Ss{t_{1}}})(T_2) = A_s + (q^3-q)^{s}(q+1)S_s + \cI(t_1, \ldots, t_s)$ was accomplished in Proposition \ref{prop:non-generic}. Observe that this element lies in the submodule $\cV \subseteq \overline{\cW}$ generated by $\cW$ and the skyscraper generators $T_{\tau(\lambda^{\pm 1}_1 \cdots \lambda^{\pm 1}_s)}$ for $t_i = \lambda_i + \lambda_i^{-1}$. This submodule $\cV$ is finitely generated and invariant under the maps $\cZg{}(L)$ and $\cZg{}(L_{[J_+]})$.

Hence, using the computations of $\cZg{}(L)$ and $\cZg{}(L_{[J_+]})$ of \cite{GP:2018a}, we can write their matricial expression in the natural set of generators of $\cV$. Now, we decompose them as $\cZg{}(L)=PDP^{-1}$ and $\cZg{}(L_{[J_+]}) = P_+D_+P_+^{-1}$, with $D$ and $D_+$ their Jordan forms (that, indeed, are diagonal matrices). Let $T_2^*=\cZg{}(D^\dag): \cV \to \Ko{\Var{\CC}}$ be the dual form of $T_2$ with respect to the standard generators of $\cV$. Then, we have that
\begin{align*}
	(q^3-q)^{g+r+s}\left[\Rep{\SL{2}(\CC)}(\Sigma_g, Q)\right] =\, & T_2^*PD^gP^{-1}P_+D_+^rP_+^{-1}\left(A_s + (q^3-q)^{s}(q+1)S_s\right) \\
	&+ T_2^*PD^gP^{-1}P_+D_+^rP_+^{-1}\left(\cI(t_1, \ldots, t_s)\right).
\end{align*}
Since $D$ and $D_+$ are diagonal matrices, the powers $D^g$ and $D_+^r$ can be computed explicitly. Hence, with the aid of a computer algebra system, the first summand gives the first term in the statement and the second summand gives the interaction term $\overline{\cI}_r(t_1, \ldots, t_s)=T_2^*PD^gP^{-1}P_+D_+^rP_+^{-1}\left(\cI(t_1, \ldots, t_s)\right)$.
\end{proof}
\end{thm}

\begin{rmk}
The formula above for $r>0$ is not valid for the case $r=0$, since $\cZg{}(L_{[J_+]}): \cV \to \cV$ has a non-trivial kernel and, thus, $D_+^0 \neq \Id$.
\end{rmk}

\begin{rmk}\label{rmk:extension-general-holonomies}
The previous result applies for a more general class of parabolic structures. Let $\Sigma_g$ be the closed orientable surface of genus $g$. Let us denote by $Q_{t, r_+, r_-}^{t_1, \ldots, t_s}$ the parabolic structure with $t$ points with holonomy $\set{-\Id}$, $r_+$ points with holonomy $[J_+]$, $r_-$ points with holonomy $[J_-]$ and $s$ points with respective holonomies $\Ss{t_1}, \ldots, \Ss{t_s}$. If $t + r_-$ is even, then $\Rep{}(\Sigma_g, Q_{t, r_+, r_-}^{t_1, \ldots, t_s})$ is isomorphic to $\Rep{}(\Sigma_g, Q_{0, r_+ + r_-, 0}^{t_1, \ldots, t_s})$, which is one of the representation varieties considered in Theorem \ref{thm:result-complete}. On the other hand, if $t+r_-$ is odd and $s > 0$, then $\Rep{}(\Sigma_g, Q_{t, r_+, r_-}^{t_1, \ldots, t_s})$ is isomorphic to $\Rep{}(\Sigma_g, Q_{0, r_+ + r_-, 0}^{-t_1, \ldots, t_s})$, which is also one of the cases considered above. The remaining case of having $t + r_-$ odd and $s=0$ can also be accomplished, as shown in \cite[Theorem 6.11]{GP:2018a}, but in general it gives a different virtual class.
\end{rmk}

\section{Parabolic $\SL{2}(\CC)$-character varieties}
\label{sec:character-varieties}

Let $G$ be a complex reductive algebraic group, let $X$ be a topological space with finitely generated fundamental group and let $Q$ be a parabolic structure, with associated representation variety $\Rep{G}(X,Q)$.

The variety $\Rep{G}(X,Q)$ parametrizes the set of all the parabolic representations $\pi_1(X) \to G$, but it does not take into account that some representations might be isomorphic. Hence, if we want to consider the moduli space of parabolic representations up to isomorphism, we need to quotient $\Rep{G}(X,Q)$ by the action of $G$ by conjugation. This can be done via the GIT quotient
$$
	\Char{G}(X,Q) = \Rep{G}(X,Q) \sslash G.
$$
This is the so-called $G$-\emph{character variety} of $X$ with parabolic structure $Q$. In this section, we will focus on the computation of the virtual class $[\Char{G}(X,Q)] \in \K{\Var{\CC}}$.

\subsection{Review of the theory of pseudo-quotients}
\label{sec:review-pseudoquotients}
In \cite{GP:2018b}, it is given an algorithm for the computation of the virtual classes of character varieties in terms of the virtual classes of the associated representation varieties. For the sake of completeness, here we will sketch briefly the method.

The key concept introduced for this aim is the pseudo-quotient. Given a complex algebraic variety $X$ with the action of an algebraic group $G$, a \emph{pseudo-quotient} of $X$ by $G$ is a complex algebraic variety $\bar{X}$ with a surjective $G$-invariant regular morphism $\pi: X \to \bar{X}$ such that, for any disjoint $G$-invariant closed sets $W_1, W_2 \subseteq X$, the Zariski closures of $\pi(W_1)$ and $\pi(W_2)$ in $\bar{X}$ are disjoint. It may be understood as a week version of a good quotient (in the GIT sense \cite{Newstead:1978}) in which we get rid of the condition that $\bar{X}$ is a categorical quotient of $X$.
Pseudo-quotients of an action may not be unique, as shown in \cite[Example 3.6]{GP:2018b}. However, in \cite[Corollary 3.10]{GP:2018b} it is proven that, if $\bar{X}_1, \bar{X}_2$ are two pseudo-quotients for the action of $G$ on $X$, then the associated virtual classes $[\bar{X}_1] = [\bar{X}_2]$, as elements of $\K{\Var{\CC}}$. In particular, if $X$ admits a GIT quotient $X \sslash G$, then $[\bar{X}] = [X \sslash G]$ for any pseudo-quotient $\bar{X}$ of $X$.

However, the key property of pseudo-quotients for our purposes is that they behave well with respect to stratifications. Suppose that $U \subseteq X$ is an orbitwise-closed open set. This means that $U$ is a Zariski open set such that, for any $x \in U$, the Zariski closure of the orbit of $x$ under the action of $G$ is also in $U$. In that case, we can decompose $X = Y \sqcup U$ with $Y$ a closed $G$-invariant set. Then, \cite[Theorem 4.1]{GP:2018b} shows that, if $\bar{Y}$ and $\bar{U}$ are pseudo-quotients for the action of $G$ on $Y$ and $U$ respectively, then $[\bar{X}] = [\bar{Y}] + [\bar{U}]$. Hence, the virtual class of a pseudo-quotient of $X$ can be computed by focusing on the corresponding ones for $Y$ and $U$.

Finally, another very useful property of pseudo-quotients holds. Suppose that there exists a subvariety $Y \subseteq X$ and an algebraic subgroup $H \subseteq G$ that `concentrates' the information of the quotient of $X$ by $G$. More precisely, this means that $H$ acts on $Y$ making it an orbitwise-closed set, that the closure of any $G$-orbit of $x \in X$ intersects with $H$ and that, for any closed $H$-invariant sets $W_1, W_2 \subseteq Y$, they are disjoint if and only if the closure of their $G$-orbits on $X$ are disjoint. This situation $(Y, H)$ is called a \emph{core} for the action. Then, in \cite[Proposition 4.4]{GP:2018b}, it is proven that if $\bar{X}$ is a pseudo-quotient of $X$ by $G$ and $\bar{Y}$ is a pseudo-quotient of $Y$ by $H$, then $[\bar{X}] = [\bar{Y}]$. Hence, we can also focus on the action on the core and, up to virtual class, it has all the information of the global quotient.

\subsection{Stratification analysis of character varieties}
With the theory of pseudo-quotients at hand, we can easily stratify the action of $G$ on the representation variety. Let us denote by $\Repred{G}(X, Q)$ and $\Repirred{G}(X, Q)$ the subvarieties of $\Rep{G}(X, Q)$ of reducible and irreducible representations, respectively. We have a decomposition $\Rep{G}(X, Q) = \Repred{G}(X, Q) \sqcup \Repirred{G}(X, Q)$ with the first stratum being a closed set and the second stratum being an open orbitwise-closed set (\cite[Proposition 6.2]{GP:2018b}). Hence, we have that
$$
	\left[\Char{G}(X, Q)\right] = \left[\Repred{G}(X, Q) \sslash G\right] + \left[\Repirred{G}(X, Q) \sslash G\right].
$$

Now, suppose that $G$ is a linear algebraic group (so in particular it is affine). Then, by \cite[Proposition 6.4]{GP:2018b} we have that $\Repirred{G}(X, Q) \to \Repirred{G}(X, Q) \sslash \mathrm{Inn}(G)$ is a free geometric quotient, where $\mathrm{Inn}(G)=G/Z(G)$ is the group of inner automorphisms of $G$. Applying now \cite[Theorem 5.4]{GP:2018b} we finally get $\left[\Repirred{G}(X, Q) \sslash G\right]\left[\mathrm{Inn}(G)\right] = \left[\Repirred{G}(X, Q)\right]$.

Furthermore, suppose that $G = \SL{n}(\CC)$. If $\Repdiag{\SL{n}(\CC)}(X, Q)$ is the subvariety of completely reduced representations with respect to the standard basis of $\CC^n$ (i.e.\ diagonal matrices), then an adaptation of Proposition 7.3 and Corollary 7.4 of \cite{GP:2018b} implies that $(\Repdiag{\SL{n}(\CC)}(X, Q), S_n)$ is a core for the action of $\SL{n}(\CC)$ on $\Repred{\SL{n}(\CC)}(X, Q)$. Here $S_n$ denotes the symmetric group that acts on $\Repdiag{\SL{n}(\CC)}(X, Q)$ by permutation of eigenvalues. Thus, $[\Repred{\SL{n}(\CC)}(X, Q) \sslash \SL{n}(\CC)] = [\Repdiag{\SL{n}(\CC)}(X) / S_n]$.

Hence, if we localize $\K{\CVar}$ by $[\PGL{n}(\CC)]=\left[\mathrm{Inn}\left(\SL{n}(\CC)\right)\right]$, these decompositions give the formula
\begin{align}\label{eqn:form-quotient-git}
	\left[\Char{\SL{n}(\CC)}(X, Q)\right] &= \left[\Repred{\SL{n}(\CC)}(X, Q) \sslash \SL{n}(\CC)\right] + \frac{\left[\Repirred{\SL{n}(\CC)}(X, Q)\right]}{\left[\PGL{n}(\CC)\right]} \\
	&= \left[\Repdiag{\SL{n}(\CC)}(X, Q) \sslash S_n\right] + \frac{\left[\Rep{\SL{n}(\CC)}(X, Q)\right]-\left[\Repred{\SL{n}(\CC)}(X, Q)\right]}{\left[\PGL{n}(\CC)\right]} \nonumber.
\end{align}

\subsection{Semi-simple punctures}
\label{sec:semi-simple-quotient}
From now on, we will consider the case $G=\SL{2}(\CC)$. In order to lighten the notation, we shall omit the group from the notation again. Let $\Sigma_g$ be the compact orientable surface of genus $g \geq 1$. Choose $t_1, \ldots, t_s \in \Bt$ and let $Q$ be the parabolic structure with $s$ punctures with holonomies $\Ss{t_1}, \ldots, \Ss{t_s}$.



In order to understand the action of $\SL{2}(\CC)$ on the representation variety $\Rep{}(\Sigma_g, Q)$, we consider $\Reput{}(\Sigma_g, Q) \subseteq \Rep{}(\Sigma_g, Q)$ the set of upper triangular representations. Let us take $A=(A_1, B_1, \ldots, A_{g}, B_g, C_1, \ldots, C_s)$ be a tuple of upper triangular matrices with $C_i \in \Ss{t_i}$, say
$$
	A = \left(
	\begin{pmatrix}
		\mu_1 & a_1 \\
		0 & \mu_1^{-1}
	\end{pmatrix},
	\begin{pmatrix}
		\nu_1 & b_1 \\
		0 & \nu_1^{-1}
	\end{pmatrix}, \ldots,
	\begin{pmatrix}
		\mu_{g} & a_{g} \\
		0 & \mu_{g}^{-1}
	\end{pmatrix},
	\begin{pmatrix}
		\nu_{g} & b_{g} \\
		0 & \nu_{g}^{-1}
	\end{pmatrix},
	\begin{pmatrix}
		\lambda_1 & c_1 \\
		0 & \lambda_1^{-1}
	\end{pmatrix}, \ldots,
	\begin{pmatrix}
		\lambda_s & c_s \\
		0 & \lambda_s^{-1}
	\end{pmatrix}
	\right)
$$
with $\mu_i, \nu_i \in \CC^*$, $\lambda_i + \lambda_i^{-1}=t_i$ and $a_i, b_i, c_i \in \CC$. Then, we have that
$$
	\prod_{i=1}^g[A_i, B_i] \prod_{k=1}^s C_k =
	\begin{pmatrix}
	{\displaystyle\lambda_1\cdots \lambda_s} & {\displaystyle \sum_{i=1}^g \mu_i\nu_i \left[\left(\nu_i-\nu_i^{-1}\right)b_i - \left(\mu_i-\mu_i^{-1}\right)a_i\right] + \sum_{k=1}^s \left(\prod_{j \neq k}\lambda_j\right) c_k= 0}
	\\ 0 & {\displaystyle\lambda_1^{-1}\cdots \lambda_s^{-1}}
	\end{pmatrix}
$$
Therefore, $A \in \Reput{}(\Sigma_g, Q)$ if and only if the following system of equations holds
\begin{equation}\label{eq:cond-upper:parabolic-diagonal}
	\left\{\begin{matrix}
	{\displaystyle \lambda_1\cdots \lambda_s = 1}, \\
	{\displaystyle \sum_{i=1}^g \mu_i\nu_i \left[\left(\nu_i-\nu_i^{-1}\right)b_i - \left(\mu_i-\mu_i^{-1}\right)a_i\right] + \sum_{k=1}^s \lambda_k^{-1} c_k= 0}.
	\end{matrix}\right.
\end{equation}

As in Section \ref{sec:interaction-term}, in order to control the first condition, we consider the set
$$
	\Lambda = \left\{(\lambda_1, \ldots, \lambda_s) \in (\CC^*)^s\,|\, \lambda_1\cdots \lambda_s =1, \lambda_i + \lambda_i^{-1} = t_i\right\},
$$
so that $\alpha_+ = \frac{1}{2}\left|\Lambda\right|$.

\begin{itemize}
	\item For the quotient of $\Repred{}(\Sigma_g, Q)$ by $\SL{2}(\CC)$, we have that $(\Repdiag{}(\Sigma_g, Q), \ZZ_2)$ is a core for the action. Since $\Repdiag{}(\Sigma_g, Q) = (\CC^*)^{2g} \times \Lambda$ and $\ZZ_2$ acts on $\Lambda$ by $(\lambda_1,\ldots, \lambda_s) \mapsto (\lambda_1^{-1},\ldots, \lambda_s^{-1})$, we obtain
$$
	\left[\Repdiag{}(\Sigma_g, Q) / \ZZ_2\right] = \left[(\CC^*)^{2g} \times \Lambda / \ZZ_2\right] = \alpha_+(q-1)^{2g}.
$$
	\item For the calculation of $\left[\Repred{}(\Sigma_g, Q)\right]$, we stratify taking care of equations (\ref{eq:cond-upper:parabolic-diagonal}).

\begin{itemize}
	\item Let $\XDh{\Rep{}(\Sigma_g, Q)} \subseteq \Repred{}(\Sigma_g, Q)$ be the set of completely reducible representations. Any element of $\XDh{\Rep{}(\Sigma_g, Q)}$ is conjugate to one of the form
$$
	\left(
	\begin{pmatrix}
		\mu_1 & 0 \\
		0 & \mu_1^{-1}
	\end{pmatrix}, \ldots,
	\begin{pmatrix}
		\nu_g & 0 \\
		0 & \nu_g^{-1}
	\end{pmatrix},
	\begin{pmatrix}
		\lambda_1 & 0 \\
		0 & \lambda_1^{-1}
	\end{pmatrix}, \ldots,
	\begin{pmatrix}
		\lambda_s & 0 \\
		0 & \lambda_s^{-1}
	\end{pmatrix}
	\right),
$$
with $(\mu_1, \ldots, \nu_g) \in (\CC^*)^{2g}$ and $(\lambda_1, \ldots, \lambda_s) \in \Lambda$. This representation is unique up to permutation of the eigenvalues, so we have a double covering
$$
	\left(\SL{2}(\CC)/{\CC^*}\right) \times (\CC^*)^{2g} \times \Lambda \longrightarrow \XDh{\Rep{}(\Sigma_g, Q)}.
$$
Therefore, as shown in \cite[Remark 5.3]{GP:2018b}, we obtain
$$
	\left[\XDh{\Rep{}(\Sigma_g, Q)}\right] = \alpha_+(q^2+q)(q-1)^{2g}.
$$
	\item Let $\XTilde{\Rep{}(\Sigma_g, Q)} \subseteq \Repred{}(\Sigma_g, Q)$ be the set of reducible but not completely reducible representations. In this case, any element is conjugated to one of the form
$$
	\left(
	\begin{pmatrix}
		\mu_1 & a_1 \\
		0 & \mu_1^{-1}
	\end{pmatrix}, \ldots,
	\begin{pmatrix}
		\nu_g & b_{g} \\
		0 & \nu_g^{-1}
	\end{pmatrix},
	\begin{pmatrix}
		\lambda_1 & c_{1} \\
		0 & \lambda_1^{-1}
	\end{pmatrix}, \ldots,
	\begin{pmatrix}
		\lambda_s & c_{s} \\
		0 & \lambda_s^{-1}
	\end{pmatrix}
	\right),
$$
with $(\mu_1, \ldots, \nu_g) \in (\CC^*)^{2g}$, $(\lambda_1, \ldots, \lambda_s) \in \Lambda$ and $(a_1, \ldots, b_{g}, c_1, \ldots, c_s) \in \pi$ where $\pi$ is the hyperplane of $\CC^{2g+s}$ given by
$$
    \pi = \left\{\sum_{i=1}^g \mu_i\nu_i \left[\left(\nu_i-\nu_i^{-1}\right)b_i - \left(\mu_i-\mu_i^{-1}\right)a_i\right] + \sum_{i=1}^s \lambda_i^{-1}c_i= 0
\right\}.
$$
In order to compute its virtual class, we have a fibration
$$
	\CC^* \times \CC \longrightarrow \PGL{2}(\CC) \times \Omega \longrightarrow \XTilde{\Rep{}(\Sigma_g, Q)}.
$$
Here, $\Omega = \left[(\CC^*)^{2g} \times \Lambda\right] \times \left[\pi -  \ell\right]$ with $\ell$ the line spanned by $(\mu_1- \mu_1^{-1}, \ldots, \nu_g - \nu_g^{-1}, \lambda_1 - \lambda_1^{-1}, \ldots, \lambda_s - \lambda_s^{-1})$. Therefore, we have
$$
	\left[\XTilde{\Rep{}(\Sigma_g, Q)}\right] = 2\alpha_+\frac{q^3-q}{(q-1)q} (q-1)^{2g}\left(q^{2g+s-1}-q\right).
$$
\end{itemize}
Hence, putting all the computations together we get
$$
	\left[\Repred{}(\Sigma_g, Q)\right] = \alpha_+ {\left(q - 1\right)}^{2  g} {\left(q + 1\right)}\left( 2q^{2  g + s - 1} -   q\right).
$$
\end{itemize}

Now observe that, since $[\PGL{2}(\CC)]=q(q+1)(q-1)$, it is invertible in $\Ko{\Var{\CC}}$. Hence, plugging all these data into formula (\ref{eqn:form-quotient-git}) and using Theorem \ref{thm:result-complete}, we finally get that
\begin{align*}
	\left[\Char{}(\Sigma_g, Q)\right] =\, & q^{2g + s - 2}(q - 1)^{2g - 2}\left(2^{2g + s - 1}- 2^s + (q + 1)^{2g + s - 2}\right) + (q^2 - 1)^{2g - 2}(q + 1)^{s}\\
	&+ 2\alpha_+(q - 1)^{2g-1}(2q - 2q^{2g + s - 2}-1)  + \frac{\overline{\cI}_0(t_1, \ldots, t_s)}{q^3-q}.
\end{align*}


\subsection{Holonomies of general type}

Now, apart from the chosen semi-simple holonomies $t_1, \ldots, t_s \in \Bt$ over punctures $p_1, \ldots, p_s \in \Sigma_g$, let us pick new different points $q_1, \ldots, q_r \in \Sigma_g$. In this case, we consider the parabolic structure $Q$ with holonomy $\Ss{t_i}$ over $p_i$, for $1 \leq i \leq s$, and $[J_+]$ over all the points $q_1, \ldots, q_r$.

Observe that, in this case, the action of $\PGL{2}(\CC)$ on the associated representation variety $\Rep{}(\Sigma_g, Q)$ is free since $\Stab_{\SL{2}(\CC)} \left(D_\lambda\right) \cap \Stab_{\SL{2}(\CC)} \left(J_+\right) = \left\{\pm \Id\right\}$. Thus, all the orbits are isomorphic and this implies that the action is also closed. Therefore, we have that
\begin{align*}
    \left[\Char{}(\Sigma_g, Q)\right] = \frac{\left[\Rep{}(\Sigma_g, Q)\right]}{\left[\PGL{2}(\CC)\right]} =\, & q^{2g + s-2}(q - 1)^{2g + r-2}\left(2^{2g+s-1} - 2^s + (q + 1)^{2g + r+s-2} \right) \\
    &+ \frac{\overline{\cI}_r(t_1, \ldots, t_s)}{q^3-q}.
\end{align*}

With these results at hand, we are ready to finish the proof of the the main theorem of this paper. 

\begin{thm}\label{thm:result-complete-quotient}
Let us fix integers $r_+, r_-, t \geq 0$, $s >0$ and traces $t_1, \ldots, t_s \in \Bt$. Consider a parabolic structure $Q$ with $r_+$ punctures with holonomy $[J_+]$, $r_-$ punctures with holonomy $[J_-]$, $s$ punctures with holonomies $\Ss{t_1}, \ldots, \Ss{t_s}$ and $t$ punctures with holonomies $\set{-\Id}$. Set $r = r_+ + r_-$ and $\sigma = \left(-1\right)^{r_- + t}$. 
The virtual class of $\Char{\SL{2}(\CC)}(\Sigma_g, Q)$ in $\Ko{\Var{\CC}}$ is
\begin{itemize}
	\item If $r>0$, then
\begin{align*}
	\left[\Char{\SL{2}(\CC)}(\Sigma_g, Q)\right] =\, & q^{2g + s-2}(q - 1)^{2g + r-2}\left(2^{2g+s-1} - 2^s + (q + 1)^{2g + r+s-2} \right) + \frac{\overline{\cI}_r(\sigma t_1, \ldots, t_s)}{q^3-q}.
\end{align*}
	\item  If $r = 0$, then
\begin{align*}
	\left[\Char{\SL{2}(\CC)}(\Sigma_g, Q)\right] =\, & q^{2g + s - 2}(q - 1)^{2g - 2}\left(2^{2g + s - 1}- 2^s + (q + 1)^{2g + s - 2}\right) + (q^2 - 1)^{2g - 2}(q + 1)^{s} \\&+ 2\alpha_{\sigma}(q - 1)^{2g-1}(2q - 2q^{2g + s - 2}-1) + \frac{\overline{\cI}_0(\sigma t_1, \ldots, t_s)}{q^3-q}.
\end{align*}
\end{itemize}

\begin{proof}
If there are no punctures with holonomies $[J_-]$ or $\set{-\Id}$, the result follows from the computations above.
In the general case, as we showed in Remark \ref{rmk:extension-general-holonomies}, if $r_-+t$ is even (i.e.\ $\sigma = 1$) we can cancel the negative holonomies $[J_-]$ and $\set{-\Id}$ pairwise and the resulting variety is isomorphic to only having $r_+ + r_-$ punctures of type $[J_+]$.

If $r_- + t$ is odd (i.e.\ $\sigma = -1$), we need to change the traces of the semi-simple holonomies to $-t_1, t_2, \ldots, t_s$. The only effect in this case is that $\alpha_+(-t_1, t_2, \ldots, t_s)=\alpha_-(t_1, t_2, \ldots, t_s)$ so the roles of $\alpha_+$ and $\alpha_-$ are interchanged.
\end{proof}
\end{thm}

If we interpret $q = [\CC] \in \K{\Var{\CC}}$ as the image in $K$-theory of the mixed Hodge structure on the compactly supported cohomology of $\CC$ (i.e.\ the Tate structure of weight $2$), the computations of this paper also give rise to the virtual Hodge structure on character varieties. Under this point of view, the results of this paper agree with the calculations of \cite{GP:2018a, GP:2018b} in the case of punctures of Jordan type.

Moreover, if we interpret $q$ as a variable in the polynomial ring $\ZZ[q, q^{-1}]$, we also obtain the $E$-polynomials of the corresponding character varieties. Under this perspective, the results of this paper agree with the existing computations in the literature. Theorem \ref{thm:result-complete-quotient} agrees with \cite{LM} for at most $2$ marked points of semi-simple holonomy and genus $1$ (notice that there is a typo in \cite{LM} in the case $r=s=1$). The result also agrees with \cite{LMN, MM:2016, MM} for a single marked point and arbitrary genus of semi-simple type. If the surface has a single marked point of type $\set{-\Id}$ (the so-called twisted variety) and arbitrary many semi-simple punctures of generic type, this result agrees with \cite{Hausel-Letellier-Villegas}. 

\begin{rmk}
The virtual classes of the character varieties with $s=0$ are also known from \cite{MM, GP:2018b}. To be precise, in these papers the computation is done for the virtual image of the mixed Hodge structure. However, following the lines of this paper, the arguments can be adapted to give the same virtual classes in $\K{\Var{\CC}}$ by interpreting $q = [\CC]$.
\end{rmk}



\bibliography{bibliography.bib}{}

\providecommand{\noopsort}[1]{}
\begin{thebibliography}{10}

\bibitem{Baraglia-Hekmati:2016}
David Baraglia and Pedram Hekmati.
\newblock Arithmetic of singular character varieties and their
  {$E$}-polynomials.
\newblock {\em Proc. Lond. Math. Soc. (3)}, 114(2):293--332, 2017.

\bibitem{Benabou}
Jean B\'enabou.
\newblock Introduction to bicategories.
\newblock In {\em Reports of the {M}idwest {C}ategory {S}eminar}, pages 1--77.
  Springer, Berlin, 1967.

\bibitem{Carlsson-Rodriguez-Villegas}
Erik Carlsson and Fernando Rodriguez~Villegas.
\newblock Vertex operators and character varieties.
\newblock {\em Adv. Math.}, 330:38--60, 2018.

\bibitem{Corlette:1988}
Kevin Corlette.
\newblock Flat {$G$}-bundles with canonical metrics.
\newblock {\em J. Differential Geom.}, 28(3):361--382, 1988.

\bibitem{Diaconescu:2017}
Duiliu-Emanuel Diaconescu.
\newblock Local curves, wild character varieties, and degenerations.
\newblock {\em Preprint arXiv:1705.05707}, 2017.

\bibitem{GPLM:2017}
{\'A}ngel Gonz\'alez-Prieto, Marina Logares, and Vicente Mu\~noz.
\newblock A lax monoidal {T}opological {Q}uantum {F}ield {T}heory for
  representation varieties.
\newblock {\em Preprint arXiv:1709.05724}, 2017.

\bibitem{Gonzalez-Prieto:Thesis}
{\'A}ngel\noopsort{01} Gonz\'alez-Prieto.
\newblock Topological quantum field theories for character varieties.
\newblock {\em PhD Thesis. Universidad Complutense de Madrid}, 2018.

\bibitem{GP:2018a}
{\'A}ngel\noopsort{02} Gonz\'alez-Prieto.
\newblock Motivic theory of representation varieties via {T}opological
  {Q}uantum {F}ield {T}heories.
\newblock {\em Preprint arXiv:1810.09714v2}, 2018.

\bibitem{GP:2018b}
{\'A}ngel\noopsort{03} Gonz\'alez-Prieto.
\newblock Stratification of algebraic quotients and character varieties.
\newblock {\em Preprint arXiv:1807.08540}, 2018.

\bibitem{Hausel:2005}
Tam\'as Hausel.
\newblock Mirror symmetry and {L}anglands duality in the non-abelian {H}odge
  theory of a curve.
\newblock In {\em Geometric methods in algebra and number theory}, volume 235
  of {\em Progr. Math.}, pages 193--217. Birkh\"auser Boston, Boston, MA, 2005.

\bibitem{Hausel-Letellier-Villegas}
Tam\'as Hausel, Emmanuel Letellier, and Fernando Rodriguez-Villegas.
\newblock Arithmetic harmonic analysis on character and quiver varieties.
\newblock {\em Duke Math. J.}, 160(2):323--400, 2011.

\bibitem{Hausel-Letellier-Villegas:2013}
Tam\'{a}s Hausel, Emmanuel Letellier, and Fernando Rodriguez-Villegas.
\newblock Arithmetic harmonic analysis on character and quiver varieties {II}.
\newblock {\em Adv. Math.}, 234:85--128, 2013.

\bibitem{Hausel-Rodriguez-Villegas:2008}
Tam\'as Hausel and Fernando Rodriguez-Villegas.
\newblock Mixed {H}odge polynomials of character varieties.
\newblock {\em Invent. Math.}, 174(3):555--624, 2008.
\newblock With an appendix by Nicholas M. Katz.

\bibitem{Hitchin:1992}
Nigel Hitchin.
\newblock Hyperk\"ahler manifolds.
\newblock In {\em S\'eminaire Bourbaki : volume 1991/92, expos\'es 745-759},
  number 206 in Ast\'erisque, pages 137--166. Soci\'et\'e math\'ematique de
  France, 1992.
\newblock talk:748.

\bibitem{LM}
Marina Logares and Vicente Mu\~noz.
\newblock Hodge polynomials of the {$\rm{SL}(2,\Bbb C)$}-character variety of
  an elliptic curve with two marked points.
\newblock {\em Internat. J. Math.}, 25(14):1450125, 22, 2014.

\bibitem{LMN}
Marina Logares, Vicente Mu\~noz, and P.~E. Newstead.
\newblock Hodge polynomials of {${\rm SL}(2,\Bbb{C})$}-character varieties for
  curves of small genus.
\newblock {\em Rev. Mat. Complut.}, 26(2):635--703, 2013.

\bibitem{MM:2016}
Javier Mart\'inez and Vicente Mu\~noz.
\newblock E-polynomials of {${SL}(2,\Bbb{C})$}-character varieties of complex
  curves of genus 3.
\newblock {\em Osaka J. Math.}, 53(3):645--681, 2016.

\bibitem{MM}
Javier Mart\'inez and Vicente Mu\~noz.
\newblock E-polynomials of the {${\rm SL}(2,\Bbb C)$}-character varieties of
  surface groups.
\newblock {\em Int. Math. Res. Not. IMRN}, (3):926--961, 2016.

\bibitem{Mellit:2017}
Anton Mellit.
\newblock Poincare polynomials of character varieties, macdonald polynomials
  and affine springer fibers.
\newblock {\em Preprint arXiv:1710.04513}, 2017.

\bibitem{Mellit}
Anton Mellit.
\newblock Poincar\'e polynomials of moduli spaces of higgs bundles and
  character varieties (no punctures).
\newblock {\em Preprint arXiv:1707.04214}, 2017.

\bibitem{Mellit:2019}
Anton Mellit.
\newblock Cell decompositions of character varieties.
\newblock {\em Preprint arXiv:1905.10685}, 2019.

\bibitem{Mereb}
Martin Mereb.
\newblock On the {$E$}-polynomials of a family of {$SL_n$}-character varieties.
\newblock {\em Math. Ann.}, 363(3-4):857--892, 2015.

\bibitem{Mozgovoy:2012}
Sergey Mozgovoy.
\newblock Solutions of the motivic {ADHM} recursion formula.
\newblock {\em Int. Math. Res. Not. IMRN}, (18):4218--4244, 2012.

\bibitem{Nakamoto}
Kazunori Nakamoto.
\newblock Representation varieties and character varieties.
\newblock {\em Publ. Res. Inst. Math. Sci.}, 36(2):159--189, 2000.

\bibitem{Newstead:1978}
P.~E. Newstead.
\newblock {\em Introduction to moduli problems and orbit spaces}, volume~51 of
  {\em Tata Institute of Fundamental Research Lectures on Mathematics and
  Physics}.
\newblock Tata Institute of Fundamental Research, Bombay; by the Narosa
  Publishing House, New Delhi, 1978.

\bibitem{Schiffmann:2016}
Olivier Schiffmann.
\newblock Indecomposable vector bundles and stable {H}iggs bundles over smooth
  projective curves.
\newblock {\em Ann. of Math. (2)}, 183(1):297--362, 2016.

\bibitem{Mozgovoy-Schiffmann}
Mozgovoy Sergey and Olivier Schiffmann.
\newblock Counting higgs bundles and type a quiver bundles.
\newblock {\em Preprint arXiv:1705.04849}, 2017.

\bibitem{Simpson:parabolic}
Carlos~T. Simpson.
\newblock Harmonic bundles on noncompact curves.
\newblock {\em J. Amer. Math. Soc.}, 3(3):713--770, 1990.

\bibitem{Simpson:1992}
Carlos~T. Simpson.
\newblock Higgs bundles and local systems.
\newblock {\em Inst. Hautes \'Etudes Sci. Publ. Math.}, (75):5--95, 1992.

\bibitem{SimpsonI}
Carlos~T. Simpson.
\newblock Moduli of representations of the fundamental group of a smooth
  projective variety. {I}.
\newblock {\em Inst. Hautes \'Etudes Sci. Publ. Math.}, (79):47--129, 1994.

\bibitem{SimpsonII}
Carlos~T. Simpson.
\newblock Moduli of representations of the fundamental group of a smooth
  projective variety. {II}.
\newblock {\em Inst. Hautes \'Etudes Sci. Publ. Math.}, (80):5--79 (1995),
  1994.

\end{thebibliography}
\bibliographystyle{plain}

\end{document}